\documentclass[12pt]{amsart}

\usepackage[usenames]{color}
\usepackage{amsmath,amsthm}
\usepackage{amssymb,latexsym}
\usepackage{enumerate}
\usepackage[usenames]{color}
\usepackage{amssymb}

\usepackage[colorlinks=true,linkcolor=webgreen,filecolor=webmaroon,
		citecolor=webred]{hyperref}
\definecolor{webgreen}{rgb}{0,.5,0}
\definecolor{webbrown}{rgb}{.6,0,0} 
\definecolor{webred}{rgb}{.9,.1,0}
\definecolor{webmaroo}{rgb}{0, 0.87, 0.68}
\usepackage{float}
\usepackage{graphics,amsmath,amssymb}
\usepackage{amsfonts}
\usepackage{amsthm}
\usepackage{graphicx}
\usepackage{epsf}

\newtheorem{theorem}{Theorem}

\newtheorem{lem}[theorem]{Lemma}

\theoremstyle{definition}

\theoremstyle{remark}

\numberwithin{equation}{section}

\DeclareMathOperator{\dd}{\thinspace d}

\DeclareMathOperator{\Li}{Li}

\begin{document}
\vskip 0.2true cm 
 
\title
[Analytic implications from the p.n.t.] 
{Analytic implications\\ from the prime number theorem}

\author[Y. Cheng]{Yuanyou Cheng} 
\address{\scriptsize 
Yuanyou Cheng\hfill\break 
\indent 
\scriptsize
Dept. of Mathematics, University of Pennsylvania, Philadelphia, PA 19104 
and Dept. of Mathematics, Brandeis University, Waltham, MA 02452, USA.
\hskip 0.2true cm\hfill\break
\indent{\bf Current Address}: 
Department of Mathematics, The University of Maryland, College 
Park, MD 20742, USA. \hfill
}
\email{\scriptsize 
Yuanyou Furui Cheng <Cheng.Y.Fred@math.umd.edu> }

\author[G. Fox]{Glenn J. Fox}
\address{\scriptsize
{\tt Glenn J. Fox},
Department of Mathematics and Physical Sciences, 
Rogers State University, Claremore, OK 74017, USA.
}
\email{\scriptsize 
Glenn J. Fox <\hskip -.02true cm gfox@rsu.edu
\hskip -.16true cm>}

\author[M. Hassani]{Mehdi Hassani } 
\address{\scriptsize
{\tt Mehdi Hassani},
Department of Mathematics, University of Zanjan,
University Blvd., 45371-38791, Zanjan, Iran}
\email{\scriptsize 
Mehdi Hassani
<\hskip -.02true cm 
mehdi.hassani@znu.ac.ir
\hskip -.16true cm>}  

\markboth{Y. Cheng, G. Fox and M. Hassani}{right head}


\thanks{*The corresponding author would like to thank Carl B. Pomerance for 
his continuous encouragement, and his helpful comments during the writing of 
this article. He wants to thank Sergio Albeverio and Andrew M. Odlyzko for 
helpful comments. Also, we wish to thank the reviewers for their careful reading.}

\subjclass[2010]{11N05, 11M26, 32A60, 11Y35, 11A41, 11R42.\hfill} 

\date{Accepted on October 5, 2014. Final version on April ??, 2021.}
 
\dedicatory{Dedicated to Carl B. Pomerance on the occasion of 
his 75$\sp{th}$ birthday}
\vfill

\maketitle


\vskip -2true cm 
\thispagestyle{empty}

\section{Introduction} \label{sec: sec1}
\noindent
It is well known that prime numbers play a central role in number 
theory. It has been known, since Riemann's famous memoir \cite{RB1} 
in 1859, that the distribution of prime numbers can be described 
by the zero-free region of the Riemann zeta function $\zeta(s)$. 
This function is a meromorphic function of the complex variable 
$s$. It has infinitely many zeros and a unique pole at $s=1$ 
with residue $1$. Let ${\mathbb C}$ denote the set of complex 
numbers. It is customary to denote $s =\sigma +it$, with $\sigma$ 
and $t$ real, for any $s \in{\mathbb C}$. For $\sigma >1$, 
the Riemann zeta function can be defined by 
\thispagestyle{empty}  
\begin{equation*}
\zeta(s) =\sum\sb{n=1}\sp{\infty} \dfrac{1}{n\sp{s}} 
=\prod\sb{p\in{\mathbb P} } \dfrac{1}{ 1 -\tfrac{1}{p\sp{s}} },
\end{equation*}
where ${\mathbb P}$ is the set of all prime numbers, with the second 
equality above being Euler's identity. One may verify Euler's identity 
from the fundamental theorem of arithmetic, which asserts 
\begin{equation}
\label{eq: new1added}
n =\prod\sb{l=1}\sp{k} p\sb{l}\sp{a\sb{l}},
\end{equation}
for every $n\in{\mathbb N}$, with $k\in{\mathbb N}$ and $p\sb{l}
\sp{a\sb{l}} \in {\mathbb P}\sp{\mathbb N}$, where ${\mathbb N} 
=\{1,2,3,\ldots\}$ is the set of natural numbers and ${\mathbb P}
\sp{\mathbb N}$ is that of all prime powers, that is $n \in {\mathbb P}
\sp{\mathbb N}$ if and only if $n =p\sp{k}$ for some prime $p$ and 
integer $k \ge 1$. 

For $\sigma >1$, we may use the logarithmic differentiation of Euler's 
identity to obtain
\begin{equation*}
-\dfrac{\zeta'(s)}{\zeta(s)} =\sum\sb{n=1}\sp{\infty} 
	\dfrac{\Lambda(n)}{n\sp{s}}. 
\end{equation*}
Here, the Mangoldt function $\Lambda$ is an arithmetic function defined 
by
\begin{equation*}
\Lambda(n) =
\begin{cases} 
\log p, \quad & n\in{\mathbb P}\sp{\mathbb N}; \\
0, \quad & \text{otherwise}, \\
\end{cases}
\end{equation*}
where ${\mathbb P}\sp{\mathbb N}$ is defined after \eqref{eq: new1added}.

We also use the notation ${\mathbb R}\sp{+}$ for the set of all positive 
real numbers. We shall use the symbol $\epsilon\in{\mathbb R}
\sp{+}$ for an arbitrary small positive real number, not necessarily 
the same at each occurrence in a given statement. Suppose that 
$g(x)$ and $h(x)$ are complex functions of the variable $x$ and 
$f(x)$ is a non-negative real-valued function of $x$. The notation 
$g(x) \trianglelefteq h(x) + B\, f(x)$ represents the fact that 
$|g(x) -h(x)|\le B f(x)$ where $B >0$ is a constant, whenever 
$x$ is sufficiently large, or $x\ge x\sb{0}$ for some fixed 
positive number $x\sb{0}$. We use the notation $g(x) \trianglerighteq 
h(x) + B\, f(x)$ in the similar way for the inequality in the opposite 
direction.

The Riemann zeta function has real zeros at $s =-2$, $-4$, $-6$, 
$\ldots$, called trivial zeros. Non-real zeros are known as 
non-trivial zeros. It is not very difficult to show that non-trivial 
zeros of $\zeta(s)$ are located in the commonly referred to as 
the critical strip $0< \sigma <1$. Some other results in this 
direction are zero-free regions in the form of 
\begin{equation}
\label{eq: zerofr}
\sigma >1 - h(t), \quad |t| >3,
\end{equation}
where $h(t)$, with $0 <h(t) \le \tfrac{1}{2}$, is a decreasing 
function of $t$. This function $h(t)$ includes 
$\tfrac{C}{\log |t|}$, 
$\tfrac{C\log\log|t|}{\log |t|}$, 
$\tfrac{C}{\log\sp{3/4 +\epsilon} |t|}$, 
and $\tfrac{C }{\log\sp{2/3}|t|\, (\log\log |t|)\sp{1/3}}$, 
where $C$ is a positive constant, which may be different in 
each situation. One may refer to any standard literature, 
e.g., \cite{CY4}, \cite{DH1}, \cite{IA1}, \cite{TE1}, and/or \cite{FK1}. 

We adopt the notation $\in\hskip -8.18true pt\sum\sb{n\le x} f(n) $, which 
means that we use the half-maximum convention to the sum function of 
the arithmetic function $f(n)$. Therefore, $f(x\sb{0}) =\tfrac{1}{2} 
\bigl(\lim\sb{x \to x\sb{0}-} f(x) +\lim\sb{x \to x\sb{0}+} f(x) \bigr)$. 
Let $x \ge 2$. We define the $\psi$-function similar to that in the literature, 
but with this half-maximum convention. We define the $\varpi$-function with this 
convention, as well. That is,
\begin{equation}
\label{eq: psivardefi}
\psi(x) =\ \in\hskip -9true pt\sum\sb{n\le x} \Lambda(n), 
\quad \varpi(x) =\ \in\hskip -9true pt\sum\sb{n\le x} 	
\bigl(\, \Lambda(n) -1 \bigr),
\end{equation}
where $n$ runs through the set of positive integers not greater than $x$. We 
remark here that any estimate on $\psi(x)$ may be converted to an estimate 
on $\varpi(x)$, and vice versa; the latter of which is needed later on. 
We notice that 
\begin{equation}
\label{eq: added2m12}
\ \in\hskip -9true pt\sum\sb{n\le x} 1 =
\begin{cases}
x -\{ x \}, &\quad x\not\in{\mathbb N}, \\
x -\tfrac{1}{2}, &\quad x\in{\mathbb N}, \\
\end{cases}
\end{equation}
where $\{x\}$ is the fractional part of $x$. From this and 
\begin{equation*}
\ \in\hskip -9true pt\sum\sb{n\le x} \Lambda(n) =\ \in\hskip -9true pt
	\sum\sb{n\le x} ( \Lambda(n) -1 ) +\ \in\hskip -9true pt
		\sum\sb{n\le x} 1,
\end{equation*}
we see that
\begin{equation}
\label{eq: added12}
\varpi(x) + x -1 \le\psi(x) \le \varpi(x) +x. 
\end{equation}

It is well known that a zero-free region of $\zeta(s)$, in the form of 
$\sigma >1 -h(t)$ and $|t| \ge 3$,  implies the prime number theorem 
in the following equivalent $\psi$-form and $\varpi$-form:
\begin{equation}
\label{eq: psivarpnt}
\psi(x) = x + O\bigl( x\sp{1 -H(x) } \log\sp{2} x\bigr), \quad  
\varpi(x) \trianglelefteq B\, x\sp{1 -H(x) } \log\sp{2} x,
\end{equation}
with an absolute positive constant $B$, where the function $H(x)$ is 
connected to $h(t)$ in a certain way. Less known is that the converse 
is also true. Actually, Tur\'an proved in 1950 that $\pi(x) =\Li(x) 
+O\bigl( x \exp( -C\sb{x} \log\sp{\frac{1}{1+a}} x )\bigr)$, 
for an $a$ with $0 <a <1$, implies the above zero-free 
region in \eqref{eq: zerofr} of $\zeta(s)$ with $h(t) =\tfrac{C\sb{t} }
{\log\sp{a} |t|}$ for $|t|\ge C\sb{0}$ with positive constants $C\sb{x}$, 
$C\sb{t}$, and $C\sb{0}$. See \cite{IA1} and \cite{TP1}. 

Corresponding to the definition of $\varpi(x)$ in \eqref{eq: psivardefi} 
and the estimate on $\varpi(x)$ in \eqref{eq: psivarpnt}, we may instead 
study the function 
\begin{equation}
\label{eq: deficalZ}
 -\dfrac{\zeta\sp{\prime}(s)}{\zeta(s)} 	
	-\zeta(s) =\sum\sb{n=1}\sp{\infty} \dfrac{\Lambda(n) -1}{n\sp{s}},
\end{equation}
as in \cite{CPRS}; this is the reason for us to have put Theorem 
\ref{thm: mthm} in the $\varpi$-form, although the $\psi$-form is 
well-known and used in the literature. The series in 
\eqref{eq: deficalZ} is convergent when $\sigma >1$. For $\sigma >0$, 
we have 
\begin{equation*}
\zeta(s) =\frac{s}{s-1} -s\int\sb{1}\sp{\infty} 
\frac{ v -\lfloor v\rfloor}{v\sp{s+1}} \dd v, 	 
\end{equation*}
where $\lfloor v\rfloor$ is the integer part of $v$.

In this article, we use the notation ${\mathbf Z}$ for the set of 
non-trivial zeros of $\zeta(s)$, which are denoted customarily by 
$\rho =\beta +i\,\gamma$ with $0 <\beta <1$. Let $\gamma\sb{0}$ 
($\approx 0.577215$) denote the Euler-Mascheroni constant (also 
called Euler's constant). It is well known that 
\begin{equation} 
\label{eq: derilzetaexpr}
\begin{split}
-\dfrac{\zeta'(s)}{\zeta(s)} & = \dfrac{1}{s-1} -\sum\sb{\rho \in{\mathsf Z}} 
		\Bigl(\dfrac{1}{s-\rho} +\dfrac{1}{\rho}\Bigr) -\log \pi  \\
&\hskip1.2true cm 
+\dfrac{\Gamma'\bigl(\tfrac{1}{2} s +1\bigr)}{2\Gamma\bigl(\tfrac{1}{2} s 
	+1\bigr)} +1 +\dfrac{\gamma\sb{0}}{2}  -\log 2.
\end{split}
\end{equation}
where the Gamma-function $\Gamma(s)$ has neither zeros nor poles 
for $\sigma >1$. Also, we remark that there is no pole at $s=1$ 
for the function of $s$ on the left hand side of \eqref{eq: deficalZ}, 
since the pole of $-\tfrac{\zeta\sp{\prime}(s)}{\zeta(s)}$ at $s =1$
and that of $\zeta(s)$ at the same point, cancel on the right hand 
side of \eqref{eq: deficalZ}. The set ${\mathbf Z}$ is the same as 
the set of poles of the function $-\tfrac{\zeta\sp{\prime}(s)}{\zeta(s)}$. 

It is known from \cite{EH1} that the Riemann zeta function for $|t| <14$ 
does not have any non-trivial zeros. From the computational perspective, we 
mention here that Xavier Gourdon uses an optimization of the Odlyzko and 
Sch\"onhage algorithm in \cite{OS1} and verifies in \cite{GX1} that 
the $10\sp{13}$ first non-trivial zeros of the Riemann zeta function 
are all simple and located on $\sigma =\tfrac{1}{2}$. Let $N(T)$ denote 
the number of non-trivial zeros of the Riemann zeta function in the region 
$0 <\sigma <1$ and $0 \le t\le T$. We use this resul1efined as above, if $|\gamma| 
<2445999556029$, as in \cite{CAGGP}. Henceforth, we let 
$T\sb{0} =2445999556027$ in this article, where the difference of $2$ 
is used for convenience with other related issues in the articles 
\cite{CAGGP}, \cite{CW1}, and \cite{CPRS}. 

We may adopt the notation $a\,{\mathbb N} +b$ for $a$ and $b\in{\mathbb Z}$, 
to denote the subset $\{ a\, n +b:\ n\in{\mathbb N}\}$ of ${\mathbb Z}$. Also, 
we denote 
\begin{equation}
\label{eq: new2added}
{\mathbb N}\sb{7} :=\{0, 1, 2, 3, 4, 5, 6\} \subseteq {\mathbb N},
\end{equation} 
in this article. Let $w \ge 5$. For $j \in{\mathbb N}\sb{7}$, we let $X\sb{0} 
=28.99$ and 
\begin{equation}
\label{eq: Hjxdefi}
H\sb{j}(x) =
\begin{cases}
\tfrac{1}{2}, \quad& 1\le x <X\sb{0}, \\
\tfrac{2}{\log\sp{(7 -j)/12} x},	\quad& x \ge X\sb{0}, \\
\end{cases}
\end{equation}
and  
\begin{equation}
\label{eq: hjtdefi}
h\sb{j}(t) =
\begin{cases}
\tfrac{1}{2}, \quad& |t| <T\sb{0}, \\
\tfrac{1}{2\, t\sp{(7 -j)w/12}}, \quad& |t| \ge T\sb{0}. \\
\end{cases}
\end{equation}

It is easy to see that both $H\sb{j}(x)$ and $h\sb{j}(t)$ are two-piece 
piece-wise differentiable functions for each $j\in {\mathbb N}\sb{7}$. 
They are monotonically decreasing and tend to $0$, respectively, as 
$x$ and $t$ tends to $\infty$, respectively for each fixed $j \in{\mathbb N}
\sb{7}$. The function $H\sb{j}(x)$ is continuous for all $x \ge 1$ and 
the function $h\sb{j}(t)$ has a unique jump discontinuity at its change
in definition at $t =T\sb{0}$, for all $j\in{\mathbb N}\sb{7}$.  
We remark that from the definitions in \eqref{eq: Hjxdefi} and 
\eqref{eq: hjtdefi} we have 
\begin{equation*}
\label{eq: Hjxhjt} 
\begin{split}
\tfrac{2}{\,\log\sp{(7 -j)/12} x} \le &H\sb{j}(x) \le\tfrac{1}{2}, 
	\quad \text{for } x >X\sb{0},  \\
\tfrac{1}{2\, U\, t\sp{(7 -j)w/12}} \le &h\sb{j}(t) \le \tfrac{1}{2}, 
	\quad\text{for } t >T\sb{0}, \\
\end{split}
\end{equation*}
for all $j\in{\mathbb N}$; otherwise, we have both $H\sb{j}(x) =\tfrac{1}{2}$
and $h\sb{j}(t) =\tfrac{1}{2}$ by their definitions in \eqref{eq: Hjxdefi} 
and \eqref{eq: hjtdefi}. The domain of $H\sb{j}(x)$ is $x\in[1, \infty)$ 
and that of $h\sb{j}(t)$ is $t\in(-\infty, \infty)$, even though we may
only use the function for $t \in[0, \infty)$ conveniently with the symmetry
property of the Riemann zeta function by the Schwarz principle.

Concerning the application of this setup of $H\sb{j}(x)$ and $h\sb{j}(t)$ 
later on, we let $f(x) =\log\sp{(j +5)/12}x -\tfrac{e(j +5)}{12} 
\log\log x$ for $x \in(e,\infty)$ here. We have $f\sp{\prime}(x) 
=\tfrac{j +5}{12 x \log\sp{(7 -j)/12} x} \bigl( 1 -\tfrac{e}
{ \log\sp{(j +5)/12} x} \bigr) =0$ for the unique critical point 
$x\sb{0} =e\sp{e\sp{12/(j+5)}}$. Note for $x =x\sb{0}$, we have 
$f(x) =0$. Checking the sign of $f\sp{\prime}(x)$, we see that 
$f(x) \ge 0$ for all $x \ge e$. Hence, 
\begin{equation}
\label{eq: s5temp}
\log\sp{(j +5)/12} x \ge \tfrac{e(j +5)}{12} \log\log x, \quad 
	x\in(e, \infty). 
\end{equation}
From the last inequality, here we remark that
\begin{equation}
\label{eq: s5tempA}
\dfrac{\log\sp{(j+5)/12} x}{\log x} =\dfrac{1}{\log\sp{(7 -j)/12} x} 
	\ge \dfrac{e(j +5)\log\log x}{12\log x},
\end{equation}
with the last expression corresponding to that used in the definition of $H\sb{j}(x)$. 

Our main result in this article is as follows. It is needed in \cite{CPRS} 
when we study the Riemann hypothesis, which states that the real parts of 
all non-trivial zeros of the Riemann zeta function are equal to $\tfrac{1}{2}$. 

\begin{theorem}
\label{thm: mthm}
Let $j\in{\mathbb N}\sb{7}$. If  
\begin{equation}
\label{eq: estvarpi}
\varpi(x) \trianglelefteq \, D\,x\sp{1 -H\sb{j}(x)} \log\sp{2} x, \quad x \ge X\sb{0},
\end{equation}
with $D =9$, then $\zeta(s)$ does not vanish when $\sigma >1- h\sb{j}(t)$, 
with $w \ge (j -3)\sp{2} +4$.
\end{theorem}

We sketch the proof of Theorem \ref{thm: mthm} in the next section, and provide 
the details afterwards. We remark here that the result in Theorem \ref{thm: mthm} 
translates estimates on the remainder term of the prime number theorem as an 
algebraic object into analytic descriptions for the zero-free region of 
the Riemann zeta function, and it plays a pivotal role in \cite{CPRS} with 
a proof of the Riemann hypothesis based on the results in other article 
mentioned in the sequence \cite{CY4}, \cite{CG1}, \cite{CY3}, \cite{CC1}, 
\cite{CAPG}, \cite{CMathematica}, \cite{CAGP}, \cite{CAGGP}, \cite{Lh0}, 
\cite{CW1}, \cite{CP0}, \cite{CFH}, this article, and \cite{CPRS}; with 
\cite{Ri} as a summary of insight and \cite{CLJ} as a sketch of the 
technical tools.

\medskip \section{Proof of Theorem \ref{thm: mthm}} \label{sec: proofs}
\noindent 
We shall show Theorem \ref{thm: mthm} by contradiction in this section. 
Assume to the contrary that there is a non-trivial zero $\rho\sp{\prime} 
=\beta\sp{\prime} +i\gamma\sp{\prime} \in{\mathbf Z}$ such that
\begin{equation}
\label{eq: contrary} 
\gamma\sp{\prime} > T\sb{0},\quad  1 -h\sb{j}(\gamma\sp{\prime}) 
	<\beta\sp{\prime} <1.  
\end{equation}
We let
\begin{equation}
\label{eq: delta0}
\beta\sb{0} =\tfrac{1}{2 (T\sb{0} )\sp{13/12}} <\tfrac{1}{5\times 10\sp{13}},
\end{equation}
and have 
\begin{equation}
\label{eq: beta0}
\beta\sp{\prime} >1 -h\sb{j}(T\sb{0} ) \ge 1 -h\sb{6}( T\sb{0} ) 
	=1 -\beta\sb{0},
\end{equation}
from \eqref{eq: contrary}, noting that $h\sb{j}(t)$ is a monotonously decreasing 
function of $t$ and increasing with respect to $j \in{\mathbb N}\sb{7}$, as 
$\gamma\sp{\prime} \ge T\sb{0}$, where $T\sb{0}$ is defined in the last paragraph 
before \eqref{eq: Hjxdefi}.

We prove \eqref{eq: multiplied} after some preparations in Section \ref{sec: sec3}. We 
definie ${\mathbf H}$, ${\mathbf H}\sb{0}$,  ${\mathbf H}\sb{1}$, ${\mathbf H}\sb{2}$, and 
${\mathbf H}\sb{3}$ as specific subsets of ${\mathbf Z}$, and sums on these subsets $S$, 
$S\sb{0}$, $S\sb{1}$, $S\sb{2}$, and $S\sb{3}$, respectively. We give upper bounds 
for $|S\sb{0}|$, $|S\sb{1}|$, $|S\sb{2}|$, and $|S\sb{3}|$, in Section \ref{sec: sec4}, by 
means of the inequality in \eqref{eq: multiplied}. An upper bound on $|S|$ is given in 
\eqref{eq: keysum}, making use of the relationship between $S$ and other four subset
$S\sb{0}$, $S\sb{1}$, $S\sb{2}$, and $S\sb{3}$ given in \eqref{eq: sto4}. We then take 
the advantage of Lemma \ref{lem: fromCWA} to acquire a lower bound for $S$ as in 
\eqref{eq: 227now}. By our assumption in \eqref{eq: contrary}, we obtain a lowerbound on $|S|$
that is greater than our upper bound. Therefore, we have proved Theorem \ref{thm: mthm} 
by contradiction.  

We define some independent constants whose values will be chosen with respect to each of 
the corresponding sub-intervals. First of all, we let 
\begin{equation}
\label{eq: sigma0}
\sigma\sb{0} >1,\quad  	s\sb{0} =\sigma\sb{0} +i\, \gamma\sp{\prime}, 
\end{equation}
where $\beta\sp{\prime}$ and $\gamma\sp{\prime}$ are subject to \eqref{eq: contrary}. 

Next, we let 
\begin{equation}
\label{eq: check}
1 <{\hat x} \le 1 +\tfrac{1 -2\beta\sb{0}}{2( \sigma\sb{0} -1 +\beta\sb{0} ) }, 
\end{equation}
for a horizontal restraint, which is needed in \eqref{eq: newreq}. Recalling 
\eqref{eq: beta0}, one sees that the second restriction in the above implies ${\hat x} 
\le \tfrac{\sigma\sb{0} -1/2}{\sigma\sb{0} -\beta\sp{\prime}}$. Also, we let 
\begin{equation}
\label{eq: udefi}
{\ddot u} =\tfrac{1}{\sigma\sb{0} -\beta\sp{\prime}},
\end{equation}
as a restraint for directions both horizontal and vertical directions. The values of 
these constants $\sigma\sb{0}$ and ${\hat x}$ will be chosen later on. For a vertical 
restriction, we let
\begin{equation}
\label{eq: hori}
1 ={\ddot u} (\sigma\sb{0} -\beta\sp{\prime}) <{\hat y} <T\sb{0},  \quad {\hat y} 
	\in{\mathbb N} +1.
\end{equation}
The values of the following constants $\sigma\sb{0}$, ${\hat x}$, ${\ddot u}$, and 
${\hat y}$ will be chosen later on.  We notice that 
\begin{equation}
\label{eq: newreq} 
\sigma\sb{0} -{\hat x} ( \sigma\sb{0} -\beta\sp{\prime} )  >\tfrac{1}{2},
\end{equation} 
from \eqref{eq: check}, recalling the remark afterwards.

Then, we define 
\begin{equation}
\label{eq: domain}
{\mathbf  H} =\Bigl\{ \rho\in{\mathbf Z}: |\gamma -\gamma\sp{\prime}| 
	\le {\ddot u}(\sigma\sb{0} -\beta\sp{\prime}) \quad \text{and}\quad 
	\beta >\sigma\sb{0} -{\hat x} ( \sigma\sb{0} -\beta\sp{\prime} )  \Bigr\},
\end{equation}
and $L$ to be the number of the zeros of the Riemann zeta function in 
the region ${\bf H}$. For convenience, we also let ${\mathbf H}\sb{0} 
={\mathbf Z}$ and define 
\begin{figure}[htb]
\begin{center}
\includegraphics[width=1\textwidth, angle=0]{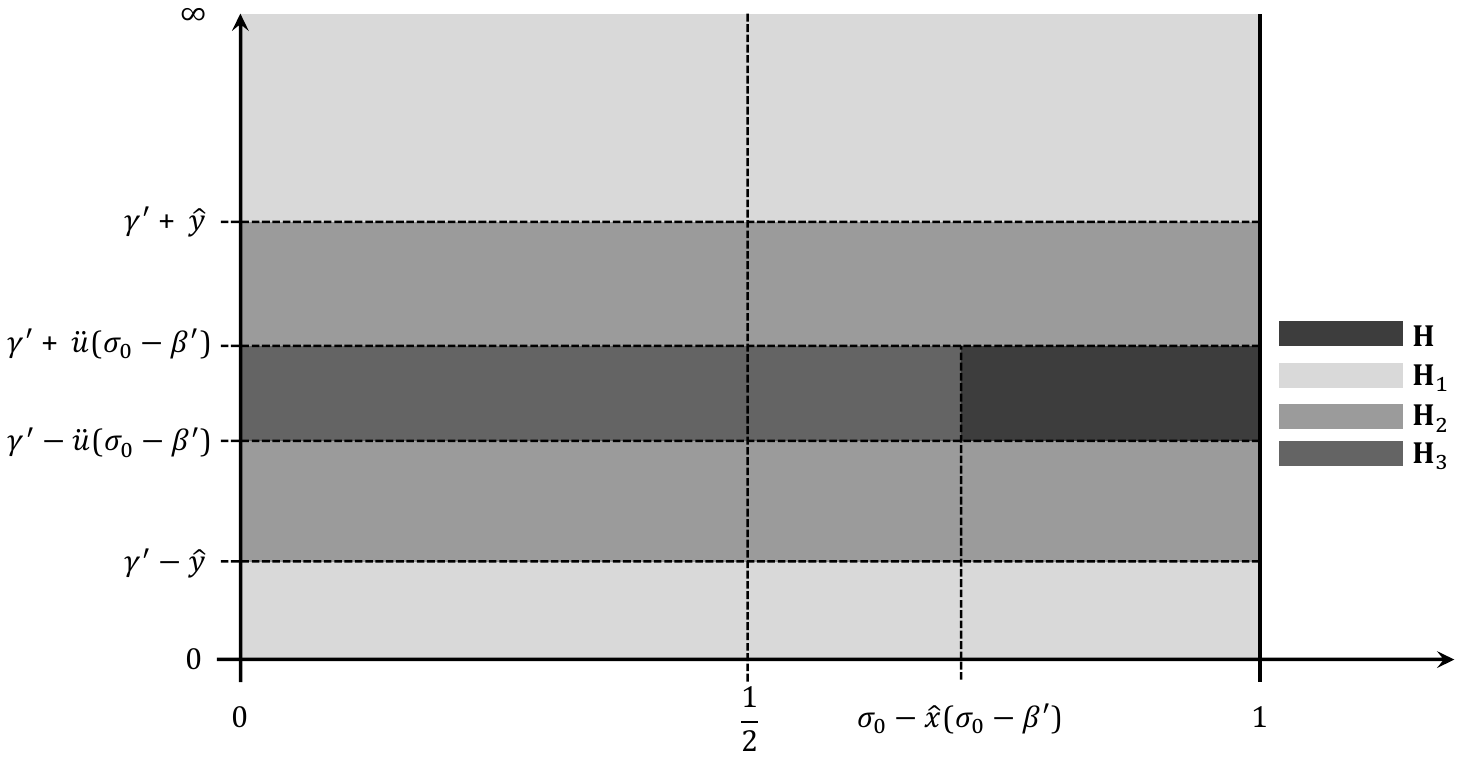}
\caption{${\mathbf H}\sb{0} ={\mathbf H}\sb{1}\cup{\mathbf H}\sb{2}
	\cup{\mathbf H}\sb{3}\cup{\mathbf H}$}
\label{fig: figlabel}
\end{center}
\end{figure}
\begin{equation}
\label{eq: threesets}
\begin{split}
&{\mathbf H}\sb{1} =\bigl\{\rho\in {\mathbf Z}:\  |\gamma - \gamma\sp{\prime}|
	>{\hat y} \bigr\},  \\
&{\mathbf H}\sb{2} =\bigl\{\rho\in {\mathbf Z}:\ {\ddot u}(\sigma\sb{0} 
	-\beta\sp{\prime}) 	<|\gamma -\gamma\sp{\prime}|\le{\hat y} \bigr\}, \\ 
&{\mathbf H}\sb{3} =\bigl\{\rho\in {\mathbf Z}:\  |\gamma -\gamma\sp{\prime}| 
	\le {\ddot u}(\sigma\sb{0} -\beta\sp{\prime})\ , 
	\beta \le\sigma\sb{0} -{\hat x} ( \sigma\sb{0} -\beta\sp{\prime} )  \bigr\}. \\
\end{split}
\end{equation} 
We remark here that ${\mathbf H}\sb{3}$ is on the left side of ${\mathbf H}$, 
whereas ${\mathbf H}$ is on the right side of the half line $\sigma =\tfrac{1}{2}$, 
which can be seen from \eqref{eq: newreq}. Also, ${\mathbf H}\sb{2}$ covers regions 
above and below the union of ${\mathbf H}\sb{3}$ and ${\mathbf H}$, and 
${\mathbf H}\sb{1}$ covers the outside regions above and below ${\mathbf H}\sb{2}$. 
In particular, we emphasize that the region ${\mathbf H}$ is completely located in 
the open half plane to the right of the line $\sigma =\tfrac{1}{2}$, 
by \eqref{eq: newreq}. We shall need this property when we give the lower bound 
of $S$ as in \eqref{eq: 227now} in Section \ref{sec: sec5}. 

Let $N(\lambda, T)$ be the number of zeros of $\zeta(s)$ when $\Re(s) \ge\lambda$ 
and $0 \le \Im(s) \le T$, with $\tfrac{1}{2} \le \lambda \le 1$ and $T \ge0$. In 
\cite{CW1}, there is a stronger estimate on the zero-growth rate equivalent to 
the Lindel\"of hypothesis, which is proved with the newly introduced pseudo-Gamma 
function by Y. Cheng, C. B. Pomerance, G. J. Fox, and S. W. Graham,  in \cite{CAPG}. 
We quote the result of Theorem 1 from \cite{CW1} as the lemma below. 

\begin{lem} 
\label{lem: fromCWA} 
Let $T \ge T\sb{0}$ with $T\sb{0}$ designed in the paragraph before 
\eqref{eq: Hjxdefi}. 
Then, for $\tfrac{1}{2} <\lambda <1$ and $1\le d \le \tfrac{5}{4}$, we have 
\begin{equation}
\label{eq: lowrate} 
N(\lambda, T+d) -N(\lambda, T-d) \le 3.
\end{equation} 
\end{lem}

Still, we use the constants $b$ and $c$, which are related to 
$\gamma\sp{\prime}$, satisfying 
\begin{equation}
\label{eq: bcrange}
b \ge\tfrac{U}{40},
\end{equation} 
in which, the restriction on $c$ is to guarantee that the value of $\tau\sb{0}$ 
in \eqref{eq: tausiota} below is positive; we actually let 
\begin{equation}
\label{eq: c}
c\log\gamma\sp{\prime} -b\log\gamma\sp{\prime} \ge U,
\end{equation}
for the application in \eqref{eq: new57}. This choice of $c$ also place a restriction 
on the choice of $b$. The values of $b$ is determined later to restrict the choice of 
a constant $k$, noting that $\gamma\sp{\prime} \ge T\sb{0}$, such that 
\begin{equation}
\label{eq: krange}
b\,\log\gamma\sp{\prime} \le k \le c\,\log\gamma\sp{\prime}, 
	\quad k \in{\mathbb N} +6. 
\end{equation}

Recall the definition of ${\mathbb N}\sb{7}$ and $j\in{\mathbb N}\sb{7}$
in \eqref{eq: new2added} so that $0\le j\le 6$. We let $J(t)$ as below and note 
\begin{equation}
\label{eq: Gadded}
J(t) =\bigl[ 6\, t\,( t +4) \bigr]\sp{\frac{6}{(j +5)e -12}}\,
	\ge\log\tfrac{11\, t\sp{2} +1}{10},
\end{equation}
as $\tfrac{6}{(j +5) e -12} \ge \tfrac{6}{11 e -12} >\tfrac{1}{3}$, $[6 t(t +4)]
\sp{\frac{6}{(j +5)e -12}} \ge t\sp{2/3}$, $t\sp{2/3} >3\log t$, and $3\log t 
>\log\tfrac{11 t\sp{2} +1}{10}$ for $t \ge T\sb{0}$, noting that the denominator 
$(j +5) e -12 \ge 5 e -12 > 1$in the exponent of the first component. The above 
two components in the maximum in \eqref{eq: Gadded} will be used for 
\eqref{eq: resultone} and \eqref{eq: sumup}, respectively. Then, we use a 
rational constant ${\mathring h}\in{\mathbb Q}$ for adjustment, such that 
\begin{equation}
\label{eq: 217}
\tfrac{1}{b k} \le {\mathring h} \le \tfrac{2}{b k}.
\end{equation}
Like other constants with convenient values that will be determined later conveniently, this 
constant ${\mathring h}$ would look superfluous once we decide on the values of $b$ and $k$. 
But, it is convenient at this point for us to use this extra constant in order to derive some 
statement more easily, especially in the sake of \eqref{eq: keykey}. 

For a technical convenience, we also let 
\begin{equation*}
\label{eq: irrational}
{\mathring\delta} ={\mathring\delta}(t)
:=\begin{cases}
0, \quad &\text{if $e\sp{ {\mathring h}\, b\, k\, J(t) } 
		\not\in{\mathbb Q}$}, \\
1, \quad &\text{if $e\sp{ {\mathring h}\, b\, k\, J(t) } 
		\in{\mathbb Q}$}. \\	
\end{cases} 
\end{equation*}
Then, we let 
\begin{equation}
\label{eq: QornotQ}	
W :=W(k; t) = e\sp{ {\mathring h}\, b\, k\, J(t) }  
	\Bigl( 1 +\dfrac{\pi}{Q} \Bigr)\sp{\mathring\delta},  
\end{equation}
where  $Q\in{\mathbb Q}$ is a sufficiently large positive 
rational number so that 
\begin{equation*}
\label{eq: Qmakeeta}
Q \ge \tfrac{\pi\,  e\sp{{\mathring h}\, b\, k\, J(t)} }{\eta},
\end{equation*} 
with respect to $\eta >0$. Furthermore, we see that $W$ is always 
an irrational number by the design of ${\mathring\delta}$, and 
\begin{equation}
\label{eq: keykey}
\begin{split}
\tfrac{ 11\, t\sp{2} +1}{10} <e\sp{t\sp{2/3} } <W 
\le( 1 +\eta)\, e\sp{2 b c \log\gamma\sp{\prime}\,
	[6\, t(t +4)]\sp{\frac{6}{e(j +5) -12} }}, \\
\end{split}
\end{equation}
for any $\eta >0$.

Now, we consider the sums $S$ and $S\sb{j}$ for $j=0$, $1$, $2$, $3$, 
where
\begin{equation}
\label{eq: sumj3}
S =\sum\sb{\rho\in{\mathbf H} }  W\sp{\rho-\rho\sp{\prime}} 
	\Bigl( \dfrac{s\sb{0} -\rho\sp{\prime}}{s\sb{0} -\rho} \Bigr)\sp{k},
\quad 
S\sb{j}=\sum\sb{\rho\in{\mathbf H}\sb{j}} W\sp{\rho-\rho\sp{\prime}} 
	\Bigl( \dfrac{s\sb{0} -\rho\sp{\prime}}{s\sb{0} -\rho} \Bigr)\sp{k},
\end{equation}
where $W =W(k; \gamma\sp{\prime})$ with $W(k; t)$ defined in \eqref{eq: QornotQ} 
above.  Moreover, ${\mathbf H}\sb{0} ={\mathbf H}\sb{1}\cup {\mathbf H}\sb{2}\cup 
{\mathbf H}\sb{3}\cup{\mathbf H}$ with	the union being disjoint, for which one may 
refer to Figure \ref{fig: figlabel} for a rough visualization. Hence,  
\begin{equation}
\label{eq: sto4}
S =S\sb{0} -S\sb{1} -S\sb{2} -S\sb{3}.
\end{equation}
We also define a constant $\tau\sp{\prime}$, which depends on $\tau\sb{0}$, 
$\tau\sb{1}$, $\tau\sb{2}$, and $\tau\sb{3}$ in \eqref{eq: tausiota} below, 
where $\tau\sb{0}$, $\tau\sb{1}$, $\tau\sb{2}$, and $\tau\sb{3}$ concerns the above 
sum $S\sb{0}$, $S\sb{1}$, $S\sb{2}$, and $S\sb{3}$, respectively, whose values are 
determined by the choices of the above independent constants later on. 

We make a preparation by giving an upper bound for a sum related to $S$ in Section 
\ref{sec: sec3}. In Section \ref{sec: sec4}, we give the estimates for $S\sb{0}$, 
$S\sb{1}$, $S\sb{2}$, and $S\sb{3}$. Recalling the definition of $\tau\sp{\prime}$ 
in \eqref{eq: tausiota}, we obtain an upper bound for the sum $S$, under the assumption 
of \eqref{eq: contrary}, in the form of  
\begin{equation}
\label{eq: keysum}
S =\sum\sb{\rho\in{\mathbf H}} \biggl( e\sp{\omega(\rho -\rho\sp{\prime})} 
	\dfrac{s\sb{0} -\rho\sp{\prime}}{s\sb{0} -\rho} \biggr)\sp{k} 
	\trianglelefteq \dfrac{ C\, W\sp{1-\beta\sp{\prime}}\, 
	\log\gamma\sp{\prime} } 
	{ \bigl(\gamma\sp{\prime}\bigr)\sp{\check\tau}}, \\
\end{equation}
where 
\begin{equation}
\label{eq: tausiota} 
{\check\tau} =\min\{\tau\sb{0}, \tau\sb{1}, \tau\sb{2}, \tau\sb{3}\}, 
\quad C =\tfrac{C\sb{0}}{28.525} +C\sb{1} +C\sb{2} +C\sb{3}, 
\end{equation}
with 
\begin{equation*}
\begin{split}
&\tau\sb{0} =1 -\tfrac{c\, \beta\sb{0} }{ \sigma\sb{0} -1 }, \quad C\sb{0} =6.7,\quad 
	\tau\sb{1} =b\log{\hat y}, \quad C\sb{1} =1.3, \\
&\tau\sb{2} =b \log{\ddot u} -\tfrac{\log[ 2.081 \delta\sb{2} ({\hat y} -1) ]}
	{\log T\sb{0}}, 	\quad C\sb{2} =2.081( 1-\delta\sb{2} ) ({\hat y} -1), \\
&\tau\sb{3} =b \log\hat x, \quad C\sb{3} =2.081, \quad 
		C =10.081+2.081(1-\delta\sb{2})({\hat y} -1), 
\end{split}
\end{equation*}
as $\gamma\sp{\prime} \ge T\sb{0}$ and $\log T\sb{0} \ge28.525$, 
where $W =W(k; t)$ with $t =\gamma\sp{\prime}$ in \eqref{eq: contrary}
and the restrictions of ${\hat y}$ in \eqref{eq: hori}, that of $
{\ddot u}$ in \eqref{eq: udefi}, and ${\check \tau}$ is defined 
in \eqref{eq: lastineq}, by $|S| \le |S\sb{0}| +|S\sb{1}| 
+|S\sb{2}| +|S\sb{3}|$ from \eqref{eq: sto4} with the upper 
bounds of $S\sb{j}$ for $j =0$, $1$, $2$, and $3$ from 
\eqref{eq: added2m12}, \eqref{eq: estS1}, \eqref{eq: estS2}, 
and \eqref{eq: estS3}. 

In Section \ref{sec: sec5}, we find a lower bound for $S$ under the assumption of 
\eqref{eq: contrary}, under the condition in \eqref{eq: bcrange}. Therefore, there exists 
at least one such sum, corresponding to a $k$ satisfying \eqref{eq: krange}, such that
\begin{equation}
\label{eq: 227now}
|S| \ge \tfrac{1}{A\sp{\prime}\,\log\sp{U} \gamma\sp{\prime}}, 
\end{equation}
where $A\sp{\prime} =\Bigl( \tfrac{ 42\, e\sp{2}\, ( b +U/\log T\sb{0} ) }{U} \Bigr)\sp{U}$, 
as $\gamma\sp{\prime}\ge\ T\sb{0}$. 

Combining \eqref{eq: keysum} and \eqref{eq: 227now}, we get 
\begin{equation}
\label{eq: finalctdn}
C\, A\sp{\prime}\, W\sp{1 -\beta\sp{\prime} }\, 	\log\sp{U}\gamma\sp{\prime} 
	\ge (\gamma\sp{\prime} )\sp{\check\tau},
\end{equation} 
or equivalently, we would have  
\begin{equation}
\label{eq: lastineq}
\begin{split}
&\hskip 1.2true cm  
-{\check\tau} \log\gamma\sp{\prime} 
	+\tfrac{ \bigl[ 6\,\gamma\sp{\prime}\, (\gamma\sp{\prime} +4) \bigr]
		\sp{ \frac{6}{e(j+5)/12 } } }{2\, (\gamma\sp{\prime})\sp{(7 -j)w/12} } \\
&\hskip 0.7true cm 
+\log \bigl[ 10.081 +2.081(1 -\delta\sb{2})({\hat y} -1) \bigr] \\
&+3\,\bigl(  \log\log\gamma\sp{\prime} +\log [ 42 e\sp{2} (b/3 +1/\log T\sb{0} ) ]  \ge 0, \\
\end{split}
\end{equation}
in which, we have used the assumption in \eqref{eq: contrary} with the definition of $h\sb{j}(t)$ 
in \eqref{eq: hjtdefi} so that $1 -\beta\sp{\prime} <h\sb{j}(\gamma\sp{\prime} ) =\tfrac{1}
{2 \,( \gamma\sp{\prime} )\sp{( 7 -j)w/2} }$ with $w \ge (j -3)\sp{2} +4$, which makes the second 
term on the left hand side of the last inequality sufficiently small and becomes a constant without 
involving $\gamma\sp{\prime}$. Hence, the above inequality is roughly equivalent to the expression
\begin{equation}
\label{eq: lastineqrough}
-{\check\tau} +0.686  +\tfrac{\log{\hat y} +\log 2}{\log T\sb{0}} \ge 0,
\end{equation}
noting
$$\dfrac{ 3\bigl[ \log\bigl(\log\gamma\sp{\prime} \times [ 42 e\sp{2} (b/3 +1/\log T\sb{0} ) ]  \bigr] }
{\log\gamma\sp{\prime}  }  \approx 0.686$$ 
and $10.081 +2.081 (1-\delta\sb{2}) ({\hat y} -1)$ 
$=2.081 ({\hat y} -1) \bigl( 1-\delta\sb{2} +\tfrac{10.081}{2.081 ({\hat y} -1) } \bigr)$ $
\approx 2 {\hat y}$, if we let $1-\delta\sb{2} +\tfrac{10.081}{2.081 ({\hat y} -1) } \approx 1$ 
by letting $\delta\sb{2} \approx\tfrac{5}{{\hat y}}$.  To finish the proof of Theorem \ref{thm: mthm},
we need to reach a contradiction by managing to have the inequality in \eqref{eq: lastineqrough}
invalid, for which we need to acquire the largest possible value for ${\check\tau}$ with suitable 
choices of all constant ``variables''. 

It is not difficult to see that we would conveniently hope to get ${\check\tau} =\tfrac{999}{1000}$ 
is only possible. It follows that $\sigma\sb{0} =1 +\tfrac{c \beta\sb{0}}{1 -\tau\sb{0}}$ by the definition 
of ${\check\tau}$ with all other involved variables known. We may tentatively choose $b =\tfrac{3}{40}$, 
which is the least value allowed for $b$. Therefore, we let $c =0.231 \ge b +\tfrac{U}{\log T\sb{0}}
\ge b +\tfrac{U}{\log \gamma\sp{\prime}}$ as required by our remark after \eqref{eq: bcrange}. 
After some experiment, we actually get ${\check\tau} =\tfrac{995}{1000}$ by setting $\sigma\sb{0} 
=1 +\tfrac{3}{2.5\times 10\sp{12}}$. 

Since the choice of ${\hat y}$ is critical for the value of $\tau\sb{1}$ and affects the value of 
$\tau\sb{2}$, we would like to choose ${\hat y}$ as small as possible in order to keep the wanted
value for ${\check\tau}$. Thus, we choose ${\hat y} =2864.073$ so that $\log{\hat y} 
\ge\tfrac{995}{1000 b} =\tfrac{199}{25}$ and $\tau\sb{1} >\tfrac{995}{1000}$ but very closely.
We recall \eqref{eq: udefi} to have ${\ddot u} =\tfrac{1}{\sigma\sb{0} -\beta\sp{\prime}} \ge 
\tfrac{1}{\sigma\sb{0} -1} =\tfrac{2.5\times 10\sp{12}}{3}$. The value of $\delta\sb{2}$ should be 
chosen as large as possible while we keep $\tau\sb{2} \ge \tfrac{995}{1000}$, for which, we 
set 
\begin{equation}
\label{eq: delta2}
\delta\sb{2} =0.143 \le \dfrac{ e\sp{b\log{\ddot u}  -995/1000}  \log T\sb{0} }{2.081 ({\hat y} -1) }. 
\end{equation}
Using of $\delta\sb{2}$ does make the final result a little better. It is easy to see that we should 
take ${\hat x} =2864.073$, as the same choice for ${\hat y}$. This choice of ${\hat x}$ also satisfies
\eqref{eq: lastineq}, as $1 +\tfrac{1 -2\beta\sb{0}}{2(\sigma\sb{0} -1+\beta\sb{0})} > 5.297\times 
10\sp{11}$.  

What remains is to compute the real final result in \eqref{eq: lastineq} with this choice of all constants 
$b$, $c$, $\sigma\sb{0}$, ${\hat x}$, ${\hat y}$, and ${\ddot u}$ to see that we really have reached 
a contradiction. In fact, we only need to check the value of the second term in \eqref{eq: lastineq} with 
\begin{equation*}
\dfrac{ \bigl[ 6\,\gamma\sp{\prime}\, (\gamma\sp{\prime} +4) \bigr]
		\sp{ \frac{6}{e(j+5)/12 } } }{2\, (\gamma\sp{\prime})\sp{(7 -j)w/12} }  
		\le \dfrac{ \bigl[ (\gamma\sp{\prime})\sp{5} \bigr]\sp{12/4.983} }{2 (\gamma\sp{\prime})\sp{13} }
		\le \dfrac{1}{\gamma\sp{\prime} } \le \dfrac{1}{T\sb{0} } <\dfrac{1}{10\sp{13}},
\end{equation*}
and 
\begin{equation*}
-\tfrac{995}{1000} +0.686 +0.3 =-0.009 <0, 
\end{equation*}
with $\tfrac{10\sp{-13}+ 10.081 +2.081 (1-\delta\sb{2}) ({\hat y} -1)}{\log T\sb{0}}< 0.3$. 
This proves that the expression on the left hand side of \eqref{eq: lastineq} is actually
negative, and this is a contradiction to \eqref{eq: lastineq} and finish the proof of 
Theorem \ref{thm: mthm}.

\section{A preparatory estimate by Tur\'an} 
\label{sec: sec3}

\noindent
In this section, we prove an inequality shown in \eqref{eq: multiplied} below
by following Tur\'an's argument in \cite{TP0}. We first prove the inequality 
\eqref{eq: result3} involving $W =W(k; t)$ below. For $s\in{\mathbb C}$ with 
$\sigma >1$ and $t \ge 25$, we have 
\begin{equation}
\label{eq: result3}
\begin{split}
&\biggl| \dfrac{W\sp{1-s}}{(s -1)\sp{k}} -\sum\sb{\rho\in{\mathbf Z}} 
\dfrac{W\sp{\rho -s}}{(s -\rho)\sp{k}} 
-\sum\sb{n=1}\sp{\infty} \dfrac{W\sp{-2n -s}}	{(s +2n)\sp{k}} \biggr| \\
&\hskip 1.6true cm  \le \tfrac{2\sp{\sigma +1} \ W\sp{1-\sigma}} {\beta\, t 
\bigl(2\sp{\sigma-1} -1\bigr)	(\sigma-1)\sp{k-1}}, \\
\end{split}
\end{equation}
under the assumption \eqref{eq: estvarpi} in Theorem \ref{thm: mthm}
stated in Section \ref{sec: sec1}. 

We recall the remark before \eqref{eq: keykey} that $W\not\in{\mathbb Q}$. 
Here, the positive integer $k$ is not less than $4$, or $k\in{\mathbb N} +3$, 
as stipulated in \eqref{eq: krange}. A similar result to \eqref{eq: result3} 
without the explicit constant, under a similar assumption 
to \eqref{eq: estvarpi}, is found in \cite{TP0}. From the inequality 
in \eqref{eq: result3}, one proves \eqref{eq: multiplied}, 
which is used in Section \ref{sec: sec4} for estimating $S$ 
in \eqref{eq: keysum}. 

In order to prove \eqref{eq: result3}, we first cite the following 
lemma from \cite{TP0}. The result in this lemma exhibits why we are 
interested in that kind of expression on the left hand 
side of \eqref{eq: result3}.
\vskip -.1true cm
\begin{lem}
\label{lem: TP0}
Let $W\in{\mathbb R}\sp{+}\backslash{\mathbb Q}$ and $k\in{\mathbb N} +3$. 
Then we have 
\begin{equation} 
\label{eq: keyeqn}
\begin{split}
&\hskip 2.2true cm 	\sum\sb{n \ge W} \dfrac{\Lambda(n)}{n\sp{s}} 
	\log\sp{k-1}\dfrac{n}{W} \\
&=(k -1 )! \Bigl[ \tfrac{W\sp{1-s}}{(s-1)\sp{k}} 
	-\sum\sb{\rho\in {\mathbf Z}} \tfrac{W\sp{\rho-s}}
	{(s-\rho)\sp{k}} -\sum\sb{n=1}\sp{\infty} 
	\tfrac{W\sp{-2n-s}}{(s+2n)\sp{k}} \Bigr],\\
\end{split} 
\end{equation} 
for $\sigma >1$.  \qed
\end{lem}

From Lemma \ref{lem: TP0}, one sees that we need to study the sum 
on the left hand side expression of \eqref{eq: keyeqn}. To do so,
we first consider the corresponding sum without ``the logarithmic factors'' 
on the left side of \eqref{eq: keyeqn}. Dividing the sum into infinitely 
many finite subsums, we have  
\begin{equation}
\label{eq: FdefiW}
F\sb{W}(s) =\sum\sb{n =N}\sp{\infty} \dfrac{\Lambda(n)}{n\sp{s}}	 
= \sum\sb{\tau=1}\sp{\infty} G\sb{\tau}(s),
\end{equation}
with 
\begin{equation}
\label{eq: for3d4}
G\sb{\tau}(s) =\sum\sb{n =N\sb{\tau}}\sp{N\sb{\tau+1} -1} 
	\dfrac{\Lambda(n)}{n\sp{s}},
\end{equation}
where $N =\lceil W\rceil$ and $N\sb{\tau} =2\sp{\tau-1}N$ 
for all $\tau\in{\mathbb N}$. For each $\tau$, one has 
\begin{equation*}
\begin{split}
&\hskip 2.5true cm G\sb{\tau}(it)  =\sum\sb{n =N\sb{\tau}}\sp{N\sb{\tau+1}-1} 
	\dfrac{\psi(n) -\psi(n-1) }{n\sp{it}} \\
&=\sum\sb{n =N\sb{\tau}}\sp{N\sb{\tau+1}-1} \psi(n) 
	\biggl( \dfrac{1}{ n\sp{it}} - \dfrac{1}{ (n+1)\sp{it}} \biggr) 
	-\biggl( \dfrac{\psi(N\sb{\tau} -1)}{N\sb{\tau}\sp{it}} 
	-\dfrac{ \psi(N\sb{\tau+1}-1) }{ N\sb{\tau+1}\sp{it}} \biggr),  \\
\end{split}
\end{equation*}
by the partial summation method. Here, we use the inequalities in 
\eqref{eq: added12} with the assumption in \eqref{eq: estvarpi}. It 
results
\begin{equation*}
\begin{split}
G\sb{\tau}(it) &\trianglelefteq \biggl| -\dfrac{N\sb{\tau} -1}
	{N\sb{\tau}\sp{it}} +\dfrac{N\sb{\tau+1}-1}
		{ N\sb{\tau+1}\sp{it}} +\sum\sb{n =N\sb{\tau}}
			\sp{N\sb{\tau+1}-1} n \biggl( \dfrac{1}
		{ n\sp{it}} -\dfrac{1}{ (n+1)\sp{it}} \biggr) \biggr| \\
&\quad +\dfrac{ \, \log\sp{2} N\sb{\tau+1} }
	{ N\sb{\tau+1}\sp{ H\sb{j}(N\sb{\tau+1} ) } } 
	\biggl( 2\, N\sb{\tau+1}  +\sum\sb{n =N\sb{\tau} }
		\sp{ N\sb{\tau+1} -1} n \biggl| \dfrac{1}{n\sp{it}} 
				-\dfrac{1}{(n+1)\sp{it}} \biggr| \biggr), \\
\end{split}
\end{equation*}
recalling the definition of $H\sb{j}(x)$ from \eqref{eq: Hjxdefi} and
noting that it is a decreasing function of $x$ for every $j\in{\mathbb N} +2$. 
The sum of the first three terms in the last expression is equal to $\sum
\sb{n =N\sb{\tau}}\sp{N\sb{\tau+1}-1} \tfrac{1}{n\sp{it}}$ by the partial 
summation method. Therefore, recalling the definition of $G\sb{\tau}(s)$
in \eqref{eq: for3d4}, we have
\begin{equation}
\label{eq: estsegment2}
\begin{split}
G\sb{\tau}(it) &\trianglelefteq 
\sum\sb{n =N\sb{\tau}}\sp{N\sb{\tau+1}-1} 
	\dfrac{1} {n\sp{it}}  +\dfrac{ \, \log\sp{2} N\sb{\tau+1} }
	{ N\sb{\tau+1}\sp{ H\sb{j}(N\sb{\tau+1} ) } } 	\\
&\quad\times \biggl( 2\, N\sb{\tau+1} +\sum\sb{n =N\sb{\tau} }
		\sp{ N\sb{\tau+1} -1} n \biggl| \dfrac{1}{n\sp{it}} 
				-\dfrac{1}{(n+1)\sp{it}} \biggr| \biggr). \\
\end{split}
\end{equation}

We need to estimate the two sums as above. 

We estimate the first sum in \eqref{eq: estsegment2} by tricky use of
the inequality $|e\sp{z} - 1 -z| \le |z|\sp{2}$ when $\Re(z) \le 1$
with $z =(1 - i\,t)\, \log\bigl( 1 +\tfrac{1}{n} \bigr)$. We also 
use $|\log( 1+u) -u| \le \tfrac{1}{2} u\sp{2}$ and $\log (1+u) 
\le u$ for $0 <u <1$. We have  
\begin{equation*}
\begin{split}
&\hskip .8true cm  \bigl| (n+1)\sp{1 -it} -n\sp{1 -it}  
	-(1 -it) n\sp{-it} \bigr| =n \Bigl| 
	\bigl(1+\tfrac{1}{n} \bigr)\sp{1 -it}  -1 -\tfrac{1 -it}{n} \Bigr| \\
&\le n\Bigl| \bigl(1+\tfrac{1}{n} \bigr)\sp{1 -it}  -1 	-(1 -it) 
	\log \bigl( 1+\tfrac{1}{n} \bigr)  \Bigr|  
	 +n |1 -it|  \Bigl| \log \bigl( 1+\tfrac{1}{n} \bigr) 
	 	-\tfrac{1}{n} \Bigr| \\
&\hskip 0.7true cm 
\le n\bigl( t\sp{2}+1 \bigr) \log\sp{2} \bigl( 1+\tfrac{1}{n} \bigr) 
	+\tfrac{\sqrt{t\sp{2}+1}}{2n} \le\tfrac{t\sp{2}+1}{n} 
		\bigl( 1 + \tfrac{1}{2\sqrt{t\sp{2} +1} } \bigr)
		< \tfrac{11 \,( t\sp{2} +1) }{10\, n}, 			
\end{split}
\end{equation*}
with $z =(1-it)\log(1+\tfrac{1}{n})$ and $u=\tfrac{1}{n}$, noting that 
$W >\tfrac{11(t\sp{2} +1)}{10}$ from \eqref{eq: keykey}, if only $t\ge 25$.

Summarizing the last inequality with its end terms from $n= N\sb{\tau}$ 
to $N\sb{\tau+1}-1$, we acquire an inequality involving the first sum
in \eqref{eq: estsegment2} with all other expressions estimable. That 
is,
\begin{equation*}
\biggl| N\sb{\tau+1}\sp{1 -it} - N\sb{\tau}\sp{1 -it} 	-(1 -it) 
\sum\sb{n=N\sb{\tau}}\sp{N\sb{\tau+1}-1} \dfrac{1}{n\sp{it}} \biggr| 
	<\dfrac{11\, (t\sp{2}+1)}{10} \sum\sb{n =N\sb{\tau}}\sp{N\sb{\tau +1} -1}
		\dfrac{1}{n}.
\end{equation*}
The last sum is bound from the above by $\tfrac{N\sb{\tau +1} }{ N\sb{\tau} } 
+\tfrac{1}{ N\sb{\tau} -1} $ from $\sum\sb{n =m}\sp{l -1} \tfrac{1}{n} \le \tfrac{1}{m -1} 
+\int\sb{m -1}\sp{l -1} \tfrac{\dd u}{u} =\tfrac{1}{m -1} +\log\tfrac{l -1}{m} \le 
\tfrac{1}{m -1} +\tfrac{l -m -1}{m} <\tfrac{1}{m -1} +\tfrac{l}{m}$ for $1< m < l$ with $m$, 
$l\in{\mathbb N}$. It follows that the last sum is bounded from above by $2 + \tfrac{1}{W -1}
<2.1$, recalling that $N\sb{\tau} \ge N >W$. From the above inequality, we then obtain 
\begin{equation}
\label{eq: sumup}
\biggl| \sum\sb{n=X\sb{\tau}}\sp{X\sb{\tau +1}-1} \dfrac{1}{n\sp{it}} \biggr|
\le 2.31\, \sqrt{t\sp{2}+1} +\dfrac{3\, N\sb{\tau}}{\sqrt{t\sp{2}+1}} 
<\dfrac{5.1\, N\sb{\tau}}{t}, 
\end{equation}
using $\bigl| |a| -|b| \bigr| \le | a -b| $, noting $N\sb{\tau +1} =2\, N\sb{\tau}$ 
and $\bigl| N\sb{\tau +1}\sp{ 1 -it} -N\sb{\tau}\sp{ 1- it} \bigr| \le N\sb{\tau +1} +N\sb{\tau}$,
and recalling $N\sb{\tau} \ge N >W(k; t) >\tfrac{11\, t\sp{2} +1}{10}$ with 
the designation of $N$ in \eqref{eq: for3d4}. 

As for the second sum in \eqref{eq: estsegment2}, we note 
\begin{equation*}
\tfrac{1}{n\sp{it}} -\tfrac{1}{(n+1)\sp{it}}
=\tfrac{1}{(n+1)\sp{it}} \Bigl( \bigl( 1 
	+\tfrac{1}{n} \bigr)\sp{it} -1 \Bigr) 
=\tfrac{1}{(n+1)\sp{it}} \Bigl( e\sp{it \log(1 + 1/n)} -1 \Bigr).
\end{equation*}
We apply the mean value theorem to the difference in the last expression, 
getting $e\sp{it\log(1+\frac{1}{n})} -1 =it\log\bigl(1+\tfrac{1}{n}\bigr) 
e\sp{i m}$ with $0< m\le t\log(1+\frac{1}{n})$. Note again that 
$\log (1 +u) \le u$ for $0< u <1$. Hence, 
\begin{equation*}
\bigl| \tfrac{1}{n\sp{it}} -\tfrac{1}{(n+1)\sp{it}} \bigr| \le \tfrac{t}{n}.
\end{equation*}
Therefore, 
\begin{equation}
\label{eq: onestep}
\sum\sb{n=N\sb{\tau}}\sp{N\sb{\tau+1} -1} n \biggl| \dfrac{1}{n\sp{it}} 
	-\dfrac{1}{(n+1)\sp{it}} \biggr| \le t\sum\sb{n =N\sb{\tau}}
		\sp{N\sb{\tau +1} -1} 1 \le N\sb{\tau}\, t, 
\end{equation}
as we have $N\sb{\tau}$ terms in the sum. 

Putting \eqref{eq: sumup} and \eqref{eq: onestep} in \eqref{eq: estsegment2}, 
we obtain
\begin{equation}
\label{eq: resultone}
\begin{split}
\bigl| G\sb{\tau}(it) \bigr| 
&\le \tfrac{ 5.1\, N\sb{\tau} }{t} +\tfrac{ \log\sp{2} N\sb{\tau+1} }
		{ N\sb{\tau+1}\sp{ H\sb{j}(N\sb{\tau+1} ) } }	\bigl( 2 N\sb{\tau +1} 
		+N\sb{\tau}\, t \bigr) \\
& =\tfrac{N\sb{\tau +1}}{2\, t} 	\Bigl( 5.1 +\tfrac{  \,t\, (t +4)\, 
	\log\sp{2} N\sb{\tau+1} }{ N\sb{\tau+1}\sp{ H\sb{j}(N\sb{\tau+1} ) } } \Bigr) 
		<\tfrac{2\sp{\tau +1} N}{t}, \\
\end{split}
\end{equation}
by 
\begin{equation*}
\dfrac{t\, (t +4)\, \log\sp{2} N\sb{\tau+1} } {N\sb{\tau+1}\sp{ H\sb{j}(N\sb{\tau+1} ) } } 
\le \dfrac{1}{6}, \ \text{or}\  e\sp{ \frac{e (j +5) -12}{6}\log\log N\sb{\tau +1} } 
\ge 6\, t\, (t +4), 
\end{equation*}
and $\tfrac{5.1 +1/6}{2} <4?$, recalling $N\sb{\tau +1} =2\sp{\tau}\, N$ for 
$\tau\in{\mathbb N}$ and the first inequality in \eqref{eq: hjtdefi} and $N\sb{\tau} >N$ 
for all $\tau\in{\mathbb N}$ after \eqref{eq: FdefiW} with $N >W$ and $W =W(k; t) \ge k\, 
e\sp{\bigl[ 6\, t (t +2) \bigr]\sp{\frac{6}{e (j +5) -12} } }$ from \eqref{eq: keykey}. 

With this estimate in \eqref{eq: resultone},  it then follows from 
$N\sb{\tau} =2\sp{\tau -1}\, N$ that 
\begin{equation}
\label{eq: subsumG}
G\sb{\tau}(s)  =\sum\sb{n =2\sp{\tau -1}N}\sp{2\sp{\tau}N -1} 
	\dfrac{\Lambda(n)}{n\sp{\sigma +it}} 
\trianglelefteq\dfrac{G\sb{\tau}(it)}{(2\sp{\tau -1}N)\sp{\sigma}}  
\le \dfrac{4}{\, 2\sp{(\sigma-1)(\tau -1)}\, N\sp{\sigma-1}\, t\, }.
\end{equation} 
Recalling the definition of $F\sb{W}(s)$ in \eqref{eq: FdefiW} with 
the definition of $G\sb{W}(s)$ in \eqref{eq: for3d4} 
with \eqref{eq: subsumG}, we acquire 
\begin{equation}
\label{eq: sumsub}
|F\sb{W}(s) |\le \dfrac{4}{t\, N\sp{\sigma -1}}\, 
	\sum\sb{j=1}\sp{\infty} \bigl( 2\sp{\sigma-1} \bigr)\sp{j -1} 
		=\dfrac{2\sp{\sigma +1} }
{t\, \bigl(2\sp{\sigma-1}-1\bigr)\, N\sp{\sigma -1} }.
\end{equation}

Next, we use the following lemma, which is from the context on page 161 
in \cite{TP1}, in estimating the expression on the left 
side of \eqref{eq: keyeqn}. 

\begin{lem}
\label{lem: contextTP}
Let $W\in{\mathbb R}\sp{+}\backslash{\mathbb Q}$ and $k\in{\mathbb N}+3$. 
Then,  
\begin{equation}
\label{eq: contextTP}
\sum\sb{n\ge W} \dfrac{\Lambda(n)}{n\sp{s}} \log\sp{k-1} \dfrac{n}{W} =(k-1) 
\int\sb{W}\sp{\infty} F\sb{u}(s) \dfrac{\log\sp{k-2} \tfrac{u}{W}}{u} \dd u,
\end{equation}
for $\sigma >1$. \qed
\end{lem}

From \eqref{eq: sumsub} and \eqref{eq: contextTP}, one has
\begin{equation*}
\label{eq: step2nd}
\begin{split}
&\hskip .8true cm \biggl| \sum\sb{n\ge W} \dfrac{\Lambda(n)}{n\sp{s}} 
	\log\sp{k-1} \dfrac{n}{W} \biggr|  \le \dfrac{4\,(k-1)\,2\sp{\sigma-1}}
		{t\, (2\sp{\sigma-1}-1)} \int\sb{W}\sp{\infty} \dfrac{\log\sp{k-2} 
		\tfrac{u}{W}}	{u\sp{\sigma}} \dd u \\
& = \dfrac{4\, (k-1)\, 2\sp{\sigma-1} W\sp{1-\sigma}}{ t (2\sp{\sigma-1}-1)} 
\int\sb{0}\sp{\infty} \dfrac{v\sp{k-2}}{e\sp{v(\sigma-1)}} \dd v 
=\dfrac{4\, (k-1)! \, 2\sp{\sigma-1} \ W\sp{1-\sigma}} { t 
	\bigl(2\sp{\sigma-1} -1\bigr)	(\sigma-1)\sp{k-1}}.  \\
\end{split}
\end{equation*}
By Lemma \ref{lem: TP0} with $(k -1)!$ in both the last expression and
\eqref{eq: keyeqn} canceled, one obtains \eqref{eq: result3}.  

To finish the proof of \eqref{eq: multiplied}, we use \eqref{eq: result3} 
to get the estimate as shown in \eqref{eq: multiplied} as below. 
Multiplying \eqref{eq: result3} with $s =s\sb{0} =\sigma\sb{0} 
+i\gamma\sp{\prime}$, or, $\sigma =\sigma\sb{0}$, and $t =\gamma\sp{\prime}$ 
and $W =W(k; \gamma\sp{\prime})$ defined in \eqref{eq: sigma0}, by a factor 
$W\sp{s\sb{0} -\rho\sp{\prime}} (s\sb{0} -\rho\sp{\prime})\sp{k}$ whose 
absolute value being $W\sp{\sigma\sb{0} -\beta\sp{\prime}} (\sigma\sb{0} -\beta\sp{\prime})\sp{k}$, recalling that $\rho\sp{\prime} 
=\beta\sp{\prime} +i\gamma\sp{\prime}$, we acquire  
\begin{equation}
\label{eq: multiplied}
\begin{split}
&W\sp{1-\rho\sp{\prime} } 	\bigl( \tfrac{s\sb{0} -\rho\sp{\prime}}{s\sb{0} -1} 
	\bigr)\sp{k} -\sum\sb{\rho\in{\mathbf Z}} W\sp{ \rho-\rho\sp{\prime} } 	
	\bigl( \tfrac{ s\sb{0} -\rho\sp{\prime} }{ s\sb{0} -\rho} \bigr)\sp{k}  
	-\sum\sb{n=1}\sp{\infty} W\sp{-2n-\rho\sp{\prime} } 
	\bigl( \tfrac{s\sb{0} -\rho\sp{\prime} }{s\sb{0} +2n} \bigr)\sp{k} \\
&\hskip 1true cm \trianglelefteq\tfrac{4\, (\sigma\sb{0}-1) 
	2\sp{\sigma\sb{0}-1} W\sp{1-\sigma\sb{0} }  } 
	{\, (2\sp{\sigma\sb{0}-1} -1) \gamma\sp{\prime}} 	
	\bigl( \tfrac{\sigma\sb{0} -\beta\sp{\prime} } 
	{\sigma\sb{0} -1} \bigr)\sp{c\log\gamma\sp{\prime}} 
	=\tfrac{4\sqrt{2} W\sp{1 -\beta\sp{\prime} } }
		   { (\gamma\sp{\prime})\sp{\tau\sb{0} } }\, \\
\end{split}
\end{equation}
with $\tau\sb{0} =1 -\tfrac{c\, \beta\sb{0}  }{  \sigma\sb{0} -1}$ as in \eqref{eq: tausiota}, since 
$\tfrac{ \sigma\sb{0} -\beta\sp{\prime}  }{ \sigma\sb{0} -1} \le \log\bigl( 1 
+\frac{ \beta\sb{0}  }{ \sigma\sb{0} -1 }  \bigr) \le \tfrac{ \beta\sb{0} }{ \sigma\sb{0} -1}$, 
recalling \eqref{eq: beta0}.

One may notice that the first sum on the left side of the last inequality 
runs over the set ${\mathbf Z}$ of all zeros  while that sum 
in \eqref{eq: keysum} involves only a subset ${\mathbf H}$ defined 
in \eqref{eq: domain} of zeros for the Riemann zeta function. That is 
the reason we have to deal with the set outside of ${\mathbf H}$ in 
the next section.

\section{Estimates on Sums}\label{sec: sec4}
\noindent  
In this section, we estimate $S\sb{j}$'s defined in \eqref{eq: sumj3} 
for $j=0$, $1$, $2$, $3$.

To estimate $S\sb{0}$, we use \eqref{eq: multiplied}, in which $S\sb{0}$
is the second sum on its left hand side. Recall the designations 
of $\sigma\sb{0}$ and $s\sb{0}$ in \eqref{eq: sigma0}, and, the restriction 
on $k$ in \eqref{eq: krange}. For the first sum on the left hand side 
of \eqref{eq: multiplied}, which is the major term in this regard, we get 
\begin{equation}
\label{eq: twoprp}
\bigl| W\sp{1-\rho\sp{\prime} } \bigl( \tfrac{s\sb{0} -\rho\sp{\prime}}
	{s\sb{0} -1}\bigr)\sp{k} \bigr| \le  W\sp{1-\beta\sp{\prime}}
	\bigl( \tfrac{\sigma\sb{0} -1/2}{\gamma\sp{\prime} } \bigr)\sp{k} 
		\le \tfrac{ W\sp{ 1-\beta\sp{\prime} } }
	{(\gamma\sp{\prime} )\sp{b \log\frac{ T\sb{0} }{\sigma\sb{0} -1/2 }  } }, 
\end{equation}
as $s\sb{0} -\rho\sp{\prime} =\sigma\sb{0} -\beta\sp{\prime} \le 
\sigma\sb{0} -\tfrac{1}{2}$ by $\beta\sp{\prime} >\tfrac{1}{2}$ 
in \eqref{eq: contrary}, $|s\sb{0} -1|>\gamma\sp{\prime}$, and, 
$k\ge b\log\gamma\sp{\prime}$. For the third sum on the left side 
of \eqref{eq: multiplied}, we have 
\begin{equation}
\label{eq: sumupd}
\begin{split}
&\hskip 1.7true cm
\biggl| \sum\sb{n=1}\sp{\infty} W\sp{-2n-\rho\sp{\prime} } 
	\biggl( \dfrac{s\sb{0} -\rho\sp{\prime} }
	{s\sb{0} +2n} \biggr)\sp{k} \biggr| \\
&\le \dfrac{1}{W\sp{2}} \biggl( \dfrac{\sigma\sb{0} -1/2}
	{\gamma\sp{\prime}} \biggr)\sp{2} \sum\sb{n=1}\sp{\infty} 
	\dfrac{1}{n\sp{k-2}} \le\dfrac{5}{4\, W\sp{2} } 
	\le \dfrac{5}{4\, e\sp{(\gamma\sp{\prime})\sp{2/3} } }\,, \\
\end{split}
\end{equation}
as $W\sp{ -2n -\rho\sp{\prime}} \trianglelefteq \tfrac{1}{W\sp{2}}$, 
$s\sb{0} +2n\trianglerighteq \gamma\sp{\prime}$, $s\sb{0} +2n 
\trianglerighteq n$, noting that $\sigma\sb{0} -\tfrac{1}{2} 
<\gamma\sp{\prime}$ from \eqref{eq: sigma0}, $k \ge 7$ from 
the stipulation after \eqref{eq: krange}, and $\sum\sb{n =1}\sp{\infty} 
\tfrac{1}{n\sp{k-2}} \le1 +\int\sb{1}\sp{\infty} \tfrac{\dd v}{v\sp{k -2} } 
=1 +\tfrac{1}{k -3} \le\tfrac{5}{4}$, and recalling 
\eqref{eq: keykey} for the larger lower bound of $W$.  
	
With \eqref{eq: multiplied}, \eqref{eq: twoprp}, and \eqref{eq: sumupd}, 
one sees that 
\begin{equation}
\label{eq: estS0}
|S\sb{0}| \le \tfrac{ 4\sqrt{2}\, W\sp{1 -\beta\sp{\prime} } }	
	{\,(\gamma\sp{\prime} )\sp{1 -\frac{c}{2(\sigma\sb{0} -1)} } } 
+\tfrac{ W\sp{1 -\beta\sp{\prime} } }{(\gamma\sp{\prime} )
	\sp{b \log \frac{T\sb{0}}{\sigma\sb{0} -1/2}  } }
+\tfrac{ 5\, W\sp{1 -\beta\sp{\prime} } }	{4\, e\sp{(\gamma\sp{\prime})\sp{2/3} } } 
	\le\tfrac{C\sb{0}\, W\sp{1-\beta\sp{\prime} } } { (\gamma\sp{\prime} )\sp{\tau\sb{0} } }, 
\end{equation}
where $C\sb{0} =6.7$ and $\tau\sb{0}$ is defined in \eqref{eq: tausiota}, 
noting that $4\sqrt{2} +1 +\tfrac{5 T\sb{0} }{4 e\sp{T\sb{0}\sp{2/3} } }
\le 6.7$, $b\log\frac{ T\sb{0} }{ \sigma -1/2 } >1$, and $\tau\sb{0} <1$ from its definition, 
as $\gamma\ge T\sb{0}$. 

For the estimate of $S\sb{1}$, we consider it by estimating its two sub-sums 
separately. One is over the set $\{\rho\in{\mathbf Z}: \gamma >\gamma\sp{\prime} 
+{\hat y} \}$; another over the set $\{\rho\in{\mathbf Z}: \gamma 
<\gamma\sp{\prime} -{\hat y} \}$. One may estimate both sub-sums similarly 
with the same upper bound, therefore, we give the details only 
for the first one. 

We need an estimate on the number of zeros in certain regions, for which we recall 
from \cite{RJ1} that for $T\ge 2$ 
\begin{equation}
\label{eq: Schoenfeld}
\bigl| N(T) - M(T) -\tfrac{7}{8} \bigr| \le Q(T) \,,
\end{equation}
where  $M(T) =\tfrac{T}{2\pi} \log\,\tfrac{T}{2\pi} -\tfrac{T}{2\pi}$ and  $Q(T) 
=0.137\log T +0.443\log\log T +1.588$. From this formula, for $t\ge 8 >2\pi +1$ 
we derive for $0< t\sb{1} <t\sb{2}$
\begin{equation}
\label{eq: f44}
N(t\sb{2}) -N(t\sb{1}) \le M(t\sb{2}) -M(t\sb{1}) +Q(t\sb{2}) +Q(t\sb{1}).
\end{equation}
We now only consider the special case, in which we use $t\sb{2} =t$ and $t\sb{1} 
=t -1$. Then, we apply the mean-value theorem for the function $f(x) =M(x)$ with 
the above function $M$ in the interval $[t -1, t]$, getting that 
the derivative $f\sp{\prime}(x) =\tfrac{1}{2\pi} \log\tfrac{x}{2\pi}$ and 
\begin{equation}
\label{eq: meanvalue1}
M( t\sb{2}) -M( t\sb{1}) \le \tfrac{1}{2\pi}  \log\tfrac{t}{2\pi}. 
\end{equation}
From \eqref{eq: f44} and \eqref{eq: meanvalue1}, we see that 
$N(t) -N(t-1) \le \tfrac{1}{2\pi}\log\tfrac{t}{2\pi} +0.137\log t +0.443\log\log t 
+1.588< 0.16 \log t +0.137\log t+0.443\log\log t$ $+1.296$ as $1.588 -\tfrac{1}{2\pi}
\log(2\pi) <1.296$. We have the following lemma. \par 

\begin{lem}
\label{lem: updPCC}
Let $t \ge 25$. Then,  
\begin{equation}
\label{eq: zparallel}
\sum\sb{\rho\in{\mathbf Z}:\ t-1 <\gamma -t \le t} 1 \le  1.04 \log t.
\end{equation}
\qed
\end{lem}

We now estimate $S\sb{j}$ defined in \eqref{eq: sumj3} with $H\sb{j}$
defined in \eqref{eq: threesets} for $j =1$, $2$, $3$. 

For $S\sb{1}$, we notice that the number of zeros $\zeta(s)$ for $\gamma
\sp{\prime}+n <\gamma \le \gamma\sp{\prime} +n +1$ is not greater than 
$1.04 \log\bigl( \gamma\sp{\prime} +n +1 \bigr)$, by the last lemma. For 
each such an $n \ge{\hat y}$, we note $\sigma\sb{0} -\beta\sp{\prime} \le 1$. 
Also, we have $s\sb{0} -\rho 
=\sigma\sb{0} + i\,\gamma\sp{\prime} -\beta -i\, \gamma$ so that 
$|s\sb{0} -\rho| \ge | \gamma\sp{\prime} -\gamma| \ge n$, recalling 
the definition of $s\sb{0}$ in \eqref{eq: sigma0}.  From the definition 
of $S\sb{1}$ in \eqref{eq: sumj3} with the definition of ${\mathbf H}\sb{1}$ 
in \eqref{eq: threesets}, we get  
\begin{equation}
\label{eq: 48m}
\begin{split}
&\hskip 0.25true cm
|S\sb{1}| \le 2 W\sp{1-\beta\sp{\prime}} \sum\sb{\rho\in{\mathbf Z}: 
	\ \gamma -\gamma\sp{\prime} >{\hat y} } \biggl( \dfrac{\sigma\sb{0} 
	-\beta\sp{\prime}}{|s\sb{0} -\rho|}  \biggr)\sp{k} \\
& \le 1.04 W\sp{1-\beta\sp{\prime}} 
	\sum\sb{n= {\hat y} +1}\sp{\infty} \dfrac{1}{n\sp{k} } 
		\log\bigl( \gamma\sp{\prime} +n+1\bigr), \\ 
\end{split}
\end{equation}
recalling that ${\hat y}\in{\mathbb N}$ from \eqref{eq: hori}. 

We note here that $\gamma\sp{\prime} +n +1 =\gamma\sp{\prime} 
\bigl( 1 +\tfrac{n +1}{\gamma\sp{\prime}} \bigr) \le 2\gamma\sp{\prime}$ 
if $n \le \gamma\sp{\prime} -1$ and $\log \gamma\sp{\prime} +\log 2 
=\log\gamma\sp{\prime} \bigl( 1 +\tfrac{\log 2}{\log\gamma\sp{\prime}} 
\bigr) \le 1.025 \log\gamma\sp{\prime}$ as $\gamma\sp{\prime} 
\ge T\sb{0}$, recalling the value of $T\sb{0}$ from the end of 
the paragraph before \eqref{eq: new2added}. This is used for the first
inequality in \eqref{eq: nbygammap}. For the sake of the second inequality
in \eqref{eq: nbygammap}, we have $\log\bigl( \gamma\sp{\prime} +n+1 \bigr) 
\le n\log \gamma\sp{\prime}$ if $(\gamma\sp{\prime})\sp{n-1} >\tfrac{n +1}
{\gamma\sp{\prime} } +1$ as $\gamma\sp{\prime} \ge T\sb{0}$; and if only 
$n\ge 2$. Also, by $u\sp{n} -u -n -1 >0$ for any $u >3$ with 
$u =\gamma\sb{\prime}$, we see that $\gamma\sp{\prime} +n +1 
\le (\gamma\sp{\prime})\sp{n}$. We get
\begin{equation}
\label{eq: nbygammap}
\log\bigl( \gamma\sp{\prime}+n+1 \bigr) \le
\begin{cases} 
1.025\log\gamma\sp{\prime}, &\quad
	\text{if ${\hat y} +1\le n <\gamma\sp{\prime} -1$ }, \\
n\,\log\gamma\sp{\prime}, &\quad
	\text{if ${\hat y} +1\le n <\infty$ }, \\
\end{cases}
\end{equation}
the latter in which only used in the case when $n \ge \lfloor \gamma\sp{\prime}
\rfloor$ in the following paragraph.

Recall the designation that ${\hat y}\ge 2$ from \eqref{eq: hori}. We let 
$\gamma\sp{\prime} -1 < g:= \lfloor\gamma\sp{\prime}\rfloor  \le \gamma\sp{\prime}$.    
\begin{equation}
\label{eq: newadded3}
\sum\sb{n ={\hat y} +1}\sp{\infty}
= \sum\sb{n ={\hat y} +1}\sp{g -1}
+
\sum\sb{n =g}\sp{\infty} =???.
\end{equation}

??
$\sum\sb{n ={\hat y} +1}\sp{\gamma\sp{\prime} -1} \tfrac{1}{n\sp{k}} 
\le \tfrac{1}{{\hat y}\sp{k}} +\int\sb{{\hat y}}\sp{\gamma\sp{\prime} -1} 
\tfrac{\dd v}{v\sp{k}} \le \tfrac{1}{{\hat y}\sp{k}} \bigl( 1 +\tfrac{{\hat y}}
{k -1} \bigr)$ and $\sum\sb{n =\lfloor\gamma\sp{\prime}\rfloor }\sp{\infty} 
\tfrac{1}{n\sp{k -1}} \le \tfrac{1}{(\gamma\sp{\prime} -1)\sp{k-1}} +\int\sb{\gamma\sp{\prime} -1 }\sp{\infty} 
\tfrac{\dd v}{v\sp{k -1}} \le \tfrac{1}{(\gamma\sp{\prime} -1)\sp{k -2}} 
\bigl( \gamma\sp{\prime} -1 +\tfrac{1}{k -2} \bigr) 
\le\bigl( 1 +\tfrac{1}{\gamma\sp{\prime}-1} \bigr)\sp{k -2} 
\tfrac{1}{(\gamma\sp{\prime})\sp{k -1}}$, as $\tfrac{\gamma\sp{\prime}}{\gamma\sp{\prime} -1}
=1+\tfrac{1}{\gamma\sp{\prime}}$ and $\gamma\sp{\prime} -1 +\tfrac{1}{k -2} <\gamma\sp{\prime}$.

Recalling \eqref{eq: 48m} with \eqref{eq: nbygammap} and the above remarks, 
we acquire
\begin{equation}
\label{eq: estS1}
\begin{split}
|S\sb{1}| &\le 1.066 \bigl( 1 +\tfrac{{\hat y}}{k -1} \bigr) 
	W\sp{1 -\beta\sp{\prime} }\, \tfrac{\log\gamma\sp{\prime}}
	{{\hat y}\sp{k} } \\
& +1.04 W\sp{1 -\beta\sp{\prime} } \tfrac{\log\gamma\sp{\prime}}
	{(\gamma\sp{\prime})\sp{k -1}} 
	\le \tfrac{ C\sb{1}\, W\sp{1 -\beta\sp{\prime} }\,
	\log\gamma\sp{\prime}}{(\gamma\sp{\prime})\sp{\tau\sb{1}}}, \\
\end{split}
\end{equation}
noting that $1.066 \bigl( 1 +\tfrac{{\hat y}}{k -1} \bigr) 
	+\bigl(1 +\tfrac{1}{\gamma\sp{\prime} -1} \bigr)\sp{k -2}
	\tfrac{1}{(\gamma\sp{\prime})\sp{k -1 -b\log{\hat y}} }$
is less than
\begin{equation*}
1.066 \bigl( 1 +\tfrac{{\hat y}}{b\log\gamma\sp{\prime} -1} \bigr) 
	+\bigl(1 +\tfrac{1}{\gamma\sp{\prime} -1} \bigr)
		\sp{b\log\gamma\sp{\prime}  -2} \tfrac{1}
			{(\gamma\sp{\prime})\sp{b\log\gamma\sp{\prime} 
				-1 -b\log{\hat y}} } <C\sb{3} =1.3, \\
\end{equation*}
as in \eqref{eq: krange}, and $\tau\sb{1}$ is defined after \eqref{eq: tausiota}. 

As for the estimate of $S\sb{2}$, we recall the definition of $S\sb{2}$ 
from \eqref{eq: sumj3} with the definition of ${\mathbf H}\sb{2}$ and
the restriction of $\hat y$ in \eqref{eq: hori}. We use \eqref{eq: zparallel} 
from Lemma \ref{lem: updPCC} to each of the $2$ subsets subject 
to $t -1 <\gamma \le t$ with $t =\gamma\sp{\prime} -{\hat y} +j$ 
for $j =1$, $2{\hat y}$. Noting that $\log \gamma\sp{\prime} +j 
\le \bigl( 1+ \tfrac{ 1+{\hat y}/T\sb{0} }{\log T\sb{0} } \bigr)
\log\gamma\sp{\prime}\le 1.000001 \log\gamma\sp{\prime}$ by ${\hat y} 
<T\sb{0}$ in \eqref{eq: hori}, for $j =1$, $2$, $\ldots$, ${\hat y} -{\ddot u} 
-1$, ${\hat y} +{\ddot u}(\sigma\sb{0} -\beta\sp{\prime})$, ${\hat y} 
+{\ddot u}(\sigma\sb{0} -\beta\sp{\prime})$, $\ldots$, $2{\hat y} -1$, we see 
that $|{\mathbf H}\sb{2}| \le 1.04\times 1.000001 \times 2\, \bigl( {\hat y} 
-{\ddot u} (\sigma\sb{0} -\beta\sp{\prime} ) \bigr)\, \log\gamma\sp{\prime} 
\le 2.081\, \bigl( {\hat y} -{\ddot u} (\sigma\sb{0} -\beta\sp{\prime} ) \bigr)\, 
\log\gamma\sp{\prime}$, recalling the set up and restrictions on ${\ddot u}$ 
and ${\hat y}$ in \eqref{eq: udefi} and \eqref{eq: hori}, respectively. 

Noting $s\sb{0} -\rho\sp{\prime} =\sigma\sb{0} -\beta\sp{\prime}$ again, recalling 
the definition of ${\ddot u}$ in \eqref{eq: udefi} and $s\sb{0} -\rho 
=\sigma\sb{0} +i\, \gamma\sp{\prime} -\beta -i\, \gamma \trianglerighteq 
|\gamma\sp{\prime} -\gamma| \ge {\ddot u}( \sigma\sb{0} -\beta\sp{\prime} )$, 
we acquire 
\begin{equation*}
|S\sb{2}| \le 2.081\, \bigl( y -{\ddot u} (\sigma\sb{0} -\beta\sp{\prime}) \bigr)\, 
	W\sp{1-\beta\sp{\prime}} \tfrac{\log\gamma\sp{\prime}}{ {\ddot u}\sp{k} }   	
		\le\tfrac{C\sb{2}\,W\sp{1-\beta\sp{\prime} }  \log\gamma\sp{\prime} }  
		{ (\gamma\sp{\prime} )\sp{\tau\sb{2}\sp{\prime} } }, 
\end{equation*}
where $C\sb{2}\sp{\prime} =2.081\bigl(  {\hat y} -{\ddot u} (\sigma\sb{0} 
-\beta\sp{\prime} ) \bigr)$ and $\tau\sb{2}\sp{\prime} =\tfrac{b}{2}\, \log{\ddot u}$ 
as in \eqref{eq: tausiota}. 

For the sake of sub-optimization, we let $0 <\delta\sb{2} <1$ and rewrite the last 
estimate in the form of 
\begin{equation}
\label{eq: estS2}
|S\sb{2}| \le\tfrac{ C\sb{2}\,W\sp{1-\beta\sp{\prime} } } 	
 				{ (\gamma\sp{\prime} )\sp{\tau\sb{2}} }, 
\end{equation}
where $C\sb{2} =2.081( 1 -\delta\sb{2} ) \bigl( {\hat y} -1 \bigr)$ and 
$\tau\sb{2} =\tau\sb{2}\sp{\prime} -\tfrac{ \log[ 2.081\, \delta\sb{2} ( {\hat y} -1 ) ] }
{ \log T\sb{0}  }$ for $\gamma\sp{\prime} \ge T\sb{0}$. This is the result in 
\eqref{eq: tausiota}.

We estimate $S\sb{3}$ similarly. Recall the definition of $S\sb{3}$ in \eqref{eq: sumj3} 
and noting that 
\begin{equation}
\label{eq: new4d12}
\tfrac{s\sb{0} -\rho\sp{\prime}}{s\sb{0} -\rho} =\tfrac{ \sigma\sb{0} -\beta\sp{\prime}}
{\sigma\sb{0} -\beta +i\, ( \gamma\sp{\prime} -\gamma)} \trianglelefteq \tfrac{1}{\hat x},
\end{equation}
from \eqref{eq: threesets} with $\sigma\sb{0} 
-\beta \ge {\hat x} ( \sigma\sb{0} -\beta\sp{\prime} )$ by the definition of ${\mathbf H}\sb{3}$, 
and Lemma \ref{lem: updPCC} with $t =\gamma\sp{\prime}$ and $t =\gamma\sp{\prime} 
+{\ddot u} (\sigma\sb{0} -\beta\sp{\prime})$ with the design of ${\ddot u}$ in \eqref{eq: udefi}, 
so that the number of zeros in ${\mathbf H}\sb{3}$ is not greater than 
$2.08\, \log (\gamma\sp{\prime} +1) \le 2.081\, \log\gamma\sp{\prime}$, 
as $\gamma\sp{\prime} \ge T\sb{0}$. We have 
\begin{equation}
\label{eq: estS3}
|S\sb{3}| \le\tfrac{C\sb{3}\, W\sp{1-\beta\sp{\prime} }  
	\log\gamma\sp{\prime} } { {\hat x}\sp{b\log \gamma\sp{\prime} } }, 
\end{equation}
where $C\sb{3} =2.081$ and $\tau\sb{3} =b \log {\hat x}$, as in \eqref{eq: tausiota}.

We finish the proof of \eqref{eq: keysum} by collecting \eqref{eq: estS0}, 
\eqref{eq: estS1}, \eqref{eq: estS2}, and \eqref{eq: estS3}, recalling 
the definition of ${\check\tau}$ in \eqref{eq: tausiota}.

\medskip
\section{Applying the power sum method lemmas} \label{sec: sec5}
\noindent
Tur\'an created the power sum method while investigating the Riemann zeta 
function and used this method to prove results about its zeros. Using his 
power sum method, Tur\'an proved the following lemma \ref{lem: keylem1} 
in \cite{TP1}. A similar result was used in \cite{TP1} by Tur\'an 
in proving similar results to Theorem \ref{thm: mthm} with respect to 
a different designation for the functions $H\sb{j}(x)$ and $h\sb{j}(t)$, 
but only in the sub-interval close to $1$. \par 

\begin{lem}
\label{lem: keylem1}
Let $L\in{\mathbb N} +1$ and $z\sb{1}$, $z\sb{2}$, 
$\ldots$, $z\sb{L}$ be complex numbers with
\begin{equation}
\label{eq: minimum} 
\min\sb{1\le l\le L} |z\sb{l}| \ge M.
\end{equation}
 Then for all $D\in{\mathbb R}\sp{+}$ such that 
 $D\ge 1$, 
\begin{equation}
\label{eq: firstlbd}
\max\sb{D \le \nu\le D+L}  \Bigl| z\sb{1}\sp{\nu} 
+z\sb{2}\sp{\nu} +\ldots +z\sb{L}\sp{\nu} \Bigr|
> M\sp{D} \biggl(\dfrac{ M \,L}{e(M+1)(D+L)}
\biggr)\sp{L}.  
\end{equation} 
\qed
\end{lem}

In Section \ref{sec: proofs}, we followed Tur\'an's proof by using 
Lemma \ref{lem: keylem1} in case (i). However, we need a slightly 
different version of Lemma \ref{lem: keylem1} in case (ii). This 
slightly different version is stated as Lemma \ref{lem: keylem2}. 

The next lemma is a modified version of Lemma T, from Turan's paper
\cite{TP2}, with minor improvements and explicit constants. Using Lemma
\ref{lem: keylem2} instead of Lemma \ref{lem: keylem1} in our application
in this section below, we have a lesser restriction on $D$ and an improved 
constant in the lower bound.

\begin{lem}
\label{lem: keylem2}
Let $L\in{\mathbb N} +1$, $l=1$, $2$, $\ldots$, 
$L$, and $z\sb{l}\in{\mathbb C}$ satisfy the condition
\begin{equation}
\label{eq: maximum}
\max\sb{1\le l\le L} |z\sb{l}| \ge 1.  
\end{equation}
Then for every $D\in{\mathbb R}\sp{+}$ such that $D\ge \tfrac{L}{40}$, 
\begin{equation}
\label{eq: lowerbd}
\max\sb{\nu:\, D\le \nu\le D+L} \bigl| z\sb{1}\sp{\nu} +z\sb{2}\sp{\nu} +\ldots +z\sb{L}\sp{\nu} \bigr| 
	\ge \biggl(\dfrac{L}{42 e\sp{2} (D+L)} \biggr)\sp{L}.   
\end{equation}
\end{lem}

We shall give the proof of Lemma \ref{lem: keylem2} in Section \ref{sec: sec6}. Here, we apply 
Lemma \ref{lem: keylem2} to validate the lower bound of $|S|$ in \eqref{eq: 227now} for at least 
one $k$subject to \eqref{eq: krange}. 

We let $l$ be a one-to-one map from ${\mathbf H}$ to $\{1, 2, \ldots, L\}$ 
such that $l =l(\rho)$ and denote 
\begin{equation}
\label{eq: zsbldefi}
z\sb{l} =z\sb{l(\rho)} =e\sp{\omega(\rho -\rho\sp{\prime})} 
	\dfrac{s\sb{0} -\rho\sp{\prime}}{s\sb{0} -\rho}.
\end{equation}
Recalling the definition of ${\mathbf H}$ in \eqref{eq: domain} with \eqref{eq: newreq}, we know 
that ${\mathbf H}$ is completely located on the right half plane $\sigma >\tfrac{1}{2}$. We apply 
Lemma \ref{lem: fromCWA} and have
\begin{equation}
\label{eq: lupper}
L =|{\mathbf H}| \le U,
\end{equation}
where $U$ is defined in Lemma \ref{lem: fromCWA}. We recall that $L >0$ from the definition 
of ${\mathbf H}$ with the assumption on contrary in \eqref{eq: contrary}. We notice that 
$z\sb{l(\rho\sp{\prime})} =1$ so that the condition in \eqref{eq: maximum} is satisfied. 

We remark here that the estimate in \eqref{eq: lowerbd} is better if $D$ is as small as possible, 
because $\bigl( \tfrac{L}{42 e\sp{2}(D+L)} \bigr)\sp{L}$ is a decreasing function of $D$ when $L$ 
is fixed. We would take $b$ to be the least possible,which also depends heavily on the choice of 
$b$, with the value required to be not less than $\tfrac{L}{40}$, to determine the value of $\check\tau$ 
in \eqref{eq: tausiota}. For this reason, we will not decide on the value of $b$ until we make 
a conclusion by putting everything together as in \eqref{eq: lastineq}. 

We then apply Lemma \ref{lem: keylem2}. We choose the quantities $D =b\log\gamma\sp{\prime}$,
which meets the requirement in Lemma \ref{lem: keylem2}, and $D +U \le c\log\gamma\sp{\prime}$, 
with $b$ and $c$ subject to the conditions in \eqref{eq: bcrange}. We notice that the expression on 
the right hand side of \eqref{eq: lowerbd} is a decreasing function of $L$ when $D$ is fixed. Therefore, 
we have $\bigl(\tfrac{L}{42 e\sp{2} (D +L)} \bigr)\sp{L} \ge\bigl( \tfrac{U}
{42 e\sp{2} (b\log\gamma\sp{\prime} +U} \bigr)\sp{U}$. By Lemma \ref{lem: keylem2}, we see that 
there exist a $\nu =k$ satisfying $b\log\gamma\sp{\prime} \le\nu \le b\log\gamma\sp{\prime} +U 
=c\log\gamma\sp{\prime}$ such that 
\begin{equation}
\label{eq: new57}
\sum\sb{l =1}\sp{L} z\sb{l}\sp{v} =\sum\sb{l =1}\sp{Ls} e\sp{ \omega(\rho -\rho\sp{\prime})} 
 	\biggl( \dfrac{s\sb{0} -\rho\sp{\prime}} {s\sb{0} -\rho} \biggr)\sp{\nu} \trianglerighteq 
 		\bigg( \dfrac{U}{42 e\sp{2}	 ( b\log\gamma\sp{\prime} +U ) } \biggr)\sp{U}, 
\end{equation}
by \eqref{eq: lupper} and \eqref{eq: krange}. However, the sum on the left hand side 
of \eqref{eq: new57} is the same as $S$ defined in \eqref{eq: sumj3} with our choice of $z\sb{l}$, 
where $\nu =k$, recalling the definition of ${\mathbf H}$ in \eqref{eq: domain}. Therefore, we 
have proved \eqref{eq: 227now}. 

This ends this section.

\medskip
\section{The power sum method lemmas} 
\label{sec: sec6}
\noindent
In this section, we prove Lemma \ref{lem: keylem2}.

Note that 
\begin{equation*}
\max\sb{\nu: D\le \nu\le D+N}  \bigl| z\sb{1}\sp{\nu} +z\sb{2}\sp{\nu} +\ldots +z\sb{L}\sp{\nu} \bigr| 
\ge \max\sb{v\in{\mathbb N}: D\le \nu\le D+L}  \bigl| z\sb{1}\sp{\nu} +z\sb{2}\sp{\nu} +\ldots 
	+z\sb{L}\sp{\nu} \bigr|,
\end{equation*}
and $\bigl(\tfrac{N}{D+N}\bigr)\sp{N}$ is a decreasing function with respect to $N$ for any fixed 
value of $D$. We see that the result in Lemma \ref{lem: keylem2} with $N$ being replaced by $L$ 
is actually stronger in the case that $N> L$. In the remain of this section, we use  Lemma 
\ref{lem: keylem1} to prove Lemma \ref{lem: keylem2} with $N$ being replaced by $L$.  For 
convenience, we use the notations  
\begin{equation}
\label{eq: 61}
\begin{split}
M\sb{0} &=\max\sb{1\le j\le L} |z\sb{j}|, \\
M\sb{1} &=\max\sb{D\le \nu\le D+L} |z\sb{1}\sp{v}
+z\sb{2}\sp{v} +\ldots +z\sb{L}\sp{v}|, \\
M\sb{2} &=\max\sb{D+1\le \nu\le D+L} |z\sb{1}\sp{v} 
+z\sb{2}\sp{v} +\ldots +z\sb{L}\sp{v}| \\
\end{split}
\end{equation}
from now on.

First of all, note that we may assume that $M\sb{0} =1$ without loss 
of generality. To justify this claim, we only need to apply Lemma 
\ref{lem: keylem2} with respect to the assumption that $M\sb{0} 
=1$ to the case in which $M\sb{0} >1$ and using $z\sb{j} / M\sb{0}$ 
in place of $z\sb{j}$ for $j=1$, $2$, $\ldots$, $L$. 
%
%

Secondly, we may assume that $D$ is an integer in Lemma \ref{lem: keylem2}
with $M\sb{1}$ defined in \eqref{eq: 61} being replaced by $M\sb{2}$ defined 
in \eqref{eq: 61}.  One may justify that the lemma is valid for any 
$D\in{\mathbb R}\sp{+}$ by using the integer part $\lfloor D\rfloor$ 
in place of $D$ and noting that $M\sb{1} \ge M\sb{2}$. 
%

We also may assume that $z\sb{j}$'s for $1\le j \le L$ are all distinct in 
proving Lemma \ref{lem: keylem2}. Otherwise, we justify the lemma by 
constructing an infinite sequence of the list $[z\sb{1k}, z\sb{2k}, \ldots, 
z\sb{Lk}]$ with respect to all $k\in{\mathbb N}$ such that the sequence 
converges to the list $[z\sb{1}, z\sb{2}, \ldots, z\sb{L}]$ and 
all $z\sb{jk}$'s are distinct for any fixed $k$ and use the limit 
$\lim\sb{k\to \infty} \max\sb{D+1\le \nu\le D+L} |z\sb{1k}\sp{v}+z\sb{2k}\sp{v} 
+\ldots +z\sb{Lk}\sp{v}|$, as in \cite{TP2}. 

Therefore, we only need to prove Lemma \ref{lem: keylem2} under the 
assumption that $z\sb{j}$'s are all distinct for $j=1$, $2$, $\ldots$, $L$, 
$M\sb{0} =1$, and $D$ is an integer with $M\sb{1}$ being replaced by 
$M\sb{2}$.

We choose 
\begin{equation}
\label{eq: hatU1}
{\hat U} =\tfrac{1}{4e \bigl( 1+\tfrac{D}{L} \bigr)} =\tfrac{L}{4 e(D +L)}, 
\end{equation} 
and we have 
\begin{equation}
\label{eq: hatU}
1-4e {\hat U} =1 -\tfrac{L}{D +L} =\tfrac{D}{D +L} >0.
\end{equation} 
We need a lemma in \cite{RS1} from the analytic theory of polynomials. 
 
\begin{lem}
\label{lem: BoutrouxC}
Let $w\in{\mathbb C}$ and $f(w) =\prod\sb{j=1}\sp{L} (w-z\sb{j})$ and 
$M\in{\mathbb R}\sp{+}$. Then for any prescribed ${\hat U}\in{\mathbb R}
\sp{+}$ the inequality $| f(w)| \ge {\hat U}\sp{L}$ holds outside at most 
$L$ discs $|w-z\sb{j}| \le r\sb{j}$ such that $r\sb{1} +r\sb{2} +\ldots 
+r\sb{L} \le 2e{\hat U}$. \qed
\end{lem}

By Lemma \ref{lem: BoutrouxC}, we have $|f(w)| \ge {\hat U}\sp{L}$ on the circle 
$|w| =r$ for some $r$ satisfying the above condition. From $|w-z\sb{j}| \le 2$ 
for every $j=1$, $2$, $\ldots$, $L$, we see that 
\begin{equation}
\label{eq: lambdap}
|w-z\sb{i\sb{1}}|\, |w-z\sb{i\sb{2}}| 
\cdots |w-z\sb{i\sb{\lambda}}| \ge
\bigl( \tfrac{{\hat U}}{2} \bigr)\sp{L}, 
\end{equation}
on $|w| =r$ for every choice of $\{ i\sb{1}, i\sb{2},
\ldots, i\sb{\lambda}\}$ from $\{1, 2, \ldots, L\}$. 
We rearrange the set $\{1, 2, \ldots, L\}$ so that we have two cases. 

Case (i). $1 =|z\sb{1}| \ge |z\sb{2}| \ge \ldots \ge  |z\sb{L}| > r$. 

Here, we use Lemma \ref{lem: keylem1} but with $M\sb{1}$ being 
replaced by $M\sb{2}$, which is valid when $D\in{\mathbb N}$ 
from \cite{TP2}. With $M =1-4e{\hat U}$ in Lemma \ref{lem: keylem1}, 
one gets 
\begin{equation}
\label{eq: estcase1}
\begin{split}
M\sb{2} &\ge (1-4e{\hat U})\sp{D} \biggl( \dfrac{(1-4e{\hat U}) L}
	{2e (1 -2e{\hat U}) (D +L)} \biggr)\sp{L} \\
&\ge (1-4e{\hat U})\sp{D} \biggl( \dfrac{L}{42 e(D +L)} 
	\biggr)\sp{L}. \\
\end{split}
\end{equation}

We recall \eqref{eq: hatU} and note that 
\begin{equation}
\label{eq: case1A}
\bigl( 1 - 4 e {\hat u} \bigr)\sp{D} =\biggl( \dfrac{D}{D +L} \biggr)\sp{D}
=\dfrac{1}{\Bigl[ \bigl( 1 + \tfrac{L}{D} \bigr)\sp{D/L} \Bigr]\sp{L} }
\ge \dfrac{1}{e\sp{L}}, 
\end{equation}
with $z =\tfrac{D}{L}$ by $(1 +z)\sp{1/z} \le e$ from the fact that the 
function $(1 +z)\sp{1/z}$ is monotonically decreasing for all $z\in(0, 
\infty)$ with $\lim\sb{z\to 0}(1 +z)\sp{1/z} =e$. Also we recall 
$\tfrac{D}{L} \ge\tfrac{1}{40}$ from the statement of Lemma \ref{lem: 
keylem2} have 
\begin{equation}
\label{eq: case1B}
\tfrac{1-4e {\hat U}}{1-2e {\hat U}} 
	=\dfrac{1 -\tfrac{1}{1 +D/L} }{1 -\tfrac{1}{2(1 +D/L)} }
	=\dfrac{1}{1 +\tfrac{1}{2(D /L)}} \ge \dfrac{1}{21}.
\end{equation}
Using these two inequalities \eqref{eq: case1A} and \eqref{eq: case1B}
in \eqref{eq: estcase1}, we deduce the estimate stated in the lemma.

Case (ii). $1=|z\sb{1}| \ge |z\sb{2}| \ge \ldots \ge |z\sb{l}| > r > |z\sb{l+1}| 
\ge \ldots \ge |z\sb{L}|$, where $l\in\{1, 2, \ldots, L-1\}$. Let 
\begin{equation*}
\label{eq: polynm}
P(w) =\prod\sb{j =l +1}\sp{L} (w -z\sb{j}) 
=\sum\sb{j=0}\sp{L-l} a\sb{j}\, w\sp{L -l -j}. 
\end{equation*}
For the coefficients of the polynomial $P(w)$, we have 
\begin{equation}
\label{eq: coefupd}
a\sb{j} = \sum\sb{l+1 \le k\sb{1} < k\sb{2} 
< \ldots < k\sb{j} \le L } z\sb{k\sb{1}}
z\sb{k\sb{2}} \cdots z\sb{k\sb{j}} 
\trianglelefteq\binom{L-l}{j} .
\end{equation}

Now, we need the following lemma, which is a classical result from 
the theory of Newton-interpolation, see page 48 in \cite{TP2}. 

\begin{lem}
\label{lem: newton}
Let $w\in{\mathbb C}$ and ${\mathcal C}$ be a simple closed 
curve consisting of analytic arcs on the $w$-plane and $G(w)$
a regular function outside and on ${\mathcal C}$ so that $G(w)
\to 0$ uniformly if $|w|\to \infty$. Let $l\in{\mathbb N}$, 
$w\sb{1}$, $w\sb{2}$, $\ldots$, $w\sb{l}$ be different 
points outside ${\mathcal C}$, and $g(w)$ be a polynomial
of degree $l-1$. If $g(w) =G(w)$ when $w =w\sb{j}$ for 
all $j=1$, $2$, $\ldots$, $l$, then 
\begin{equation*}
\label{eq: interpolate}
g(w) =\sum\sb{j =0}\sp{l-1} 
b\sb{j} \prod\sb{k=1}\sp{j} (w -w\sb{k}),
\end{equation*}
with the coefficients
\begin{equation*}
\label{eq: coeffs}
b\sb{j} =\dfrac{1}{2\pi i} \int\sb{\mathcal C} 
\dfrac{G(z)}{\prod\sb{k =1}\sp{j} (z -w\sb{k})} \dd z,
\end{equation*}
where the product $\prod\sb{k =1}\sp{0} (z -w\sb{k})$ is regarded 
to be $1$. \qed
\end{lem}

Let $Q(w)$ be the polynomial of degree $l-1$ such that 
$Q(z\sb{j}) =\tfrac{1}{z\sb{j}\sp{D +1} P(z\sb{j})}$
for every $j=1$, $2$, $\ldots$, $l$. Then, by Lemma
\ref{lem: newton}, we have 
\begin{equation}
\label{eq: coefQ}
Q(w) = \sum\sb{j=0}\sp{l-1} b\sb{j} \prod\sb{k=1}\sp{j}
(w-z\sb{k}) =\sum\sb{j=0}\sp{l-1} c\sb{j} w\sp{j} , 
\end{equation}
with 
\begin{equation*}
\label{eq: coefupdQ}
b\sb{j} = \dfrac{1}{2\pi i} \int\sb{|z| =r} \dfrac
{\dd z}{z\sp{D+1} P(z) \prod\sb{k=1}\sp{j}
(z -z\sb{k})} ,
\end{equation*}
for $j=0$, $1$, $\ldots$, $l-1$. 
From this, one gets
\begin{equation}
\label{eq: bjupd}
|b\sb{j}|  \le \dfrac{1}{r\sp{D}} \biggl( \dfrac{2}{{\hat U}}
\biggr)\sp{L} \le \dfrac{1}{(1-4eD)\sp{D}} \biggl( \dfrac{2}{{\hat U}}
\biggr)\sp{L},
\end{equation}
recalling \eqref{eq: lambdap}. 
Expressing $c\sb{j}$ in terms of $b\sb{j}$ in \eqref{eq: coefQ}  
by the above lemma, we see 
\begin{equation*}
\label{eq: cfromb}
\begin{split}
& c\sb{j} = b\sb{j} -b\sb{j+1} \sum\sb{1\le i\sb{1}\le j+1} 
z\sb{i\sb{1}} +b\sb{j+2} \sum\sb{1\le i\sb{1}, i\sb{2}\le j+1} 
z\sb{i\sb{1}} z\sb{i\sb{2}}  -\ldots \\
&\hskip .25true cm  +(-1)\sp{l-j-1} b\sb{l-1}
\sum\sb{1\le i\sb{1}, i\sb{2}, \ldots, i\sb{l-j-1}\le j+1} 
z\sb{i\sb{1}} z\sb{i\sb{2}} \cdots z\sb{i\sb{l-j-1}}, \\
\end{split} 
\end{equation*}
for $j=0$, $1$, $\ldots$, $l-2$ and $c\sb{l-1} =b\sb{l-1}$. 
By this inequality and \eqref{eq: bjupd}, we acquire
\begin{equation}
\label{eq: cjupd}
|c\sb{j}| \le \binom{l}{j+1}\dfrac{1}{(1-4eD)\sp{D}} 
\biggl( \dfrac{2}{{\hat U}} \biggr)\sp{L}
\end{equation}
recalling $|z\sb{j}|\le 1$ for $j=0$, $1$, $\ldots$, $l-2$ and noting 
that $1+\binom{j+1}{1}+\binom{j+2}{2} +\ldots +\binom{l-1}
{l-j-1} =\binom{l}{j+1}$. 
 
Finally we let 
\begin{equation}
\label{eq: defiR}
R(w) = w\sp{D +1} P(w) Q(w) =
\sum\sb{j =D +1}\sp{D +L} d\sb{j}\, w\sp{j}.
\end{equation}
It follows from the definition of $P(w)$ 
and $Q(w)$ and $z\sb{j} \ne z\sb{k}$ 
for $j \ne k$ that $R(z\sb{j}) =1$ for 
$j=1$, $2$, $\ldots$, $l$ and $R(z\sb{j})
=0$ for $j=l+1$, $l+2$, $\ldots$, $L$. 
Replacing $1$ by $R(z\sb{j})$ for all 
$j =1$, $2$, $\ldots$, $l$ in \eqref{eq: defiR}
and adding the results together, one gets 
\begin{equation} 
\label{eq: defiF}
M\sb{2}  \sum\sb{j =D+1}\sp{D+L} |d\sb{j}| \ge 1.   
\end{equation}
By \eqref{eq: defiR} with \eqref{eq: cjupd} and 
\eqref{eq: coefupd}, we obtain
\begin{equation*}
\label{eq: lastsum}
\sum\sb{j=D+1}\sp{D+L} |d\sb{j}|  
\le \biggl( \sum\sb{j=0}\sp{l-1} |c\sb{j}| \biggr) 
\biggl( \sum\sb{j=0}\sp{ L -l} |a\sb{j}| \biggr) 
\le \dfrac{1}{(1-4e{\hat U})\sp{D}} \biggl(\dfrac{4}{{\hat U}}
\biggr)\sp{L},
\end{equation*}
from which and \eqref{eq: defiF}, we deduce that  
\begin{equation}
\label{eq: estcase2}
M\sb{2} \ge (1 -4e{\hat U})\sp{D} \biggl(\dfrac{{\hat U}}{4}\biggr)\sp{L}. 
\end{equation}

Finally, we recall the choice of ${\hat U}$ in \eqref{eq: hatU1} and we see 
that the last expression is the same as the last expression 
in \eqref{eq: estcase1} as in Case (i). This finishes estimating 
in Case (ii).
 
This ends the proof of Lemma \ref{lem: keylem2}.  \qed

\medskip
\bigskip\bibliographystyle{amsplain}

\end{document}